\documentclass[12pt,twoside]{article}

\usepackage[reqno,tbtags]{amsmath}
\usepackage{amssymb}
\def\thebibliography#1{\vspace{0.14cm} 
 {\bf\Large References}\list
 {[\arabic{enumi}]}{\settowidth\labelwidth{[#1]}\leftmargin\labelwidth
 \advance\leftmargin\labelsep
 \usecounter{enumi}}
 \def\newblock{\hskip .11em plus .33em minus -.07em}
 \sloppy
 \sfcode`\.=1000\relax}

\protect\newcounter{theoremnumber}
\protect\newcounter{propositionnumber}
\protect\newcounter{lemmanumber}
\protect\newcounter{corollarynumber}
\protect
\protect
\protect\newcounter{examplenumber}
\newskip\halflineskip
\def\zerocounters{%
\protect\setcounter{corollarynumber}{0}%
\protect\setcounter{theoremnumber}{0}%
\protect\setcounter{lemmanumber}{0}%
\protect\setcounter{examplenumber}{0}%
\protect\setcounter{propositionnumber}{0}%
}%

\renewcommand{\thetheoremnumber}{\thesection.\arabic{theoremnumber}}

\newenvironment{theorem}%
{\par\addvspace{\halflineskip}
\refstepcounter{theoremnumber}
{\protect\setcounter{corollarynumber}{0}}
{\noindent \bf Theorem~\thetheoremnumber } %
\hskip.5em\ignorespaces\it}%
{\vskip\halflineskip\par\rm}

{\par\addvspace{\halflineskip}
\refstepcounter{theoremnumber}
\protect\setcounter{corollarynumber}{0}
{\noindent \bf Proposition~\thesection.\thetheoremnumber }%
\hskip.5em\ignorespaces\it}%
{\vskip\halflineskip\par\rm}

\newenvironment{lemma}%
{\par\addvspace{\halflineskip}
\refstepcounter{theoremnumber}
\protect\setcounter{corollarynumber}{0}
{\noindent \bf Lemma~\thetheoremnumber }%
\hskip.5em\ignorespaces\it}%
{\vskip\halflineskip\par\rm}

\newenvironment{cor}%
{\par\addvspace{\halflineskip}
\refstepcounter{theoremnumber}
{\noindent \bf Corollary~\thetheoremnumber }%
\hskip.5em\ignorespaces\it}%
{\vskip\halflineskip\par\rm}

\renewcommand{\thedefinitionnumber}{\thesection.\arabic{definitionnumber}}

{\par\addvspace{\halflineskip}
\refstepcounter{definitionnumber}
{\noindent \bf Definition~\thesection.\thedefinitionnumber }%
\hskip.5em\ignorespaces}%
{\vskip\halflineskip\par}

\newenvironment{example}%
{\par\addvspace{\halflineskip}
\refstepcounter{theoremnumber}
{\noindent \bf Example~\thetheoremnumber }%
\hskip.5em\ignorespaces}%
{\vskip\halflineskip\par}

\newenvironment{remark}%
{\par\addvspace{\halflineskip}
\refstepcounter{theoremnumber}
{\noindent \bf Remark~\thetheoremnumber }%
\hskip.5em\ignorespaces}%
{\vskip\halflineskip\par}

{\par\addvspace{\halflineskip}}%
{\vskip\halflineskip\par}

\def\@begintheorem#1#2{\par\addvspace{\halflineskip}
{\bf #1\  #2.}
\hskip.5em}
\def\@opargbegintheorem#1#2#3{\par\vskip\halflineskip
 {\bf #1\ #2.\ {\rm (#3)}.}\hskip.5em}
\def\@endtheorem{\vskip\halflineskip\par}

\newenvironment{proof}{\noindent{\bf Proof. }\rm}
{\unskip\nobreak\hfil\penalty50\hskip1em\hbox{}
\nobreak\hfill\qed\par\smallskip}
\def\qed{\vrule height1ex width1ex depth0pt}

{{\raisebox{0.4ex}{$\bigtriangledown$}}\vskip\halflineskip\par}

\protect

\newenvironment{acknowledgements}%
{\protect\nopagebreak\section*{\protect\raggedright \Large
Acknowledgements}}%
{}

{\protect\nopagebreak\section*{\protect\raggedright \Large
Acknowledgement}}% 
{}

\protect
\protect\newcommand{\sectionhead}[1]%
{\section {#1}}

\def\abstract{\vspace{-0.5cm}
{\bf \footnotesize Abstract.}\footnotesize }

\title{{\Large \bf Bosonic realizations of the color analogue of the
Heisenberg Lie algebra}}

\author{Gunnar Sigurdsson \vspace{0.3cm}\\ 
{\footnotesize \sl Mathematical Physics, Department of Physics, 
Royal Institute of Technology,}\\
{\footnotesize \sl AlbaNova University Center, SE-106 91 Stockholm, Sweden  } \\ 
{\footnotesize \sl gunnsi@theophys.kth.se }
\vspace{0.5cm}
\\ Sergei D. Silvestrov \vspace{0.3cm}\\
{\footnotesize \sl Centre for Mathematical Sciences, Department of
Mathematics,} \\ {\footnotesize \sl Lund Institute of Technology,
P.O. Box 118, SE-221 00 Lund, Sweden }\\ 
{\footnotesize \sl sergei.silvestrov@math.lth.se}}
\vspace{0.5cm}

\date{November 12, 2003}
\begin{document}

\maketitle
\footnote{Mathematics Subject Classification 2000: Primary 17B75;
Secondary 05A10,05A19, \\ 16G99,16S32,34K99,81S05,81S99  }

\begin{abstract}
We describe realizations of the color analogue of the Heisenberg Lie 
algebra by power series in non-commuting indeterminates satisfying 
Heisenberg's canonical commutation relations of quantum mechanics. 
The obtained formulas are used to construct new operator representations 
of the color analogue of the Heisenberg Lie algebra. These representations 
are shown to be closely connected with some combinatorial identities and 
functional difference-differential interpolation formulae involving Euler, 
Bernoulli and Stirling numbers.
\end{abstract}

\newpage
\thispagestyle{empty}
\ 
\newpage

\section{Introduction}
\label{sec:intro}
\markboth{Introduction}{Introduction}

The main object studied in this paper is the unital associative algebra 
with three generators $A_1$, $A_2$ and $A_3$ satisfying defining 
commutation relations 
\begin{align}
  A_1 A_2+A_2 A_1 &= A_3,  \notag \\
  A_1 A_3+A_3 A_1 &= 0,    \label{rel:colLHeis} \\
  A_2 A_3+A_3 A_2 &= 0.    \notag
\end{align} 
The main goal is to show how $A_1$, $A_2$ and $A_3$ can be expressed, 
using generators $A$, $B$ of the Heisenberg algebra, 
obeying Heisenberg's canonical commutation relation 
\begin{equation} 
  AB-BA=I.  \label{eq:Heisrel}
\end{equation}

The canonical representation of the commutation relation (\ref{eq:Heisrel}) 
is given by choosing $A$ as usual differentiation operator and $B$ as 
multiplication by $x$ acting on differentiable functions of one real variable 
$x$, on polynomials in one variable, or on some other suitable 
linear space of functions invariant under these operators.
In quantum mechanics, these operators, when considered on the Hilbert space 
of square integrable functions, are essentially the same as the canonical 
Heisenberg-Schr{\"o}dinger observables of momentum and coordinate, differing 
just by a complex scaling factor. The Heisenberg canonical commutation 
relation (\ref{eq:Heisrel}) is also satisfied by the annihilation and creation 
operators in a quantum harmonic oscillator.
Whenever $A_1$, $A_2$ and $A_3$ satisfy (\ref{rel:colLHeis}) and are 
interpreted as observables within some physical system, the problem we 
consider is that of realization of these observables within a physical system 
described by the Heisenberg-Schr{\"o}dinger observables or by the quantum 
harmonic oscillator model. 
This point of view can be very valuable for physical applications as a step 
towards understanding bosonic realizations of fermionic, super-symmetric
or color systems.
     
A complex associative algebra $L$ with generators $A_1$, $A_2$, $A_3$ and 
defining relations \mbox{(\ref{rel:colLHeis})} is called the graded analogue 
of the Heisenberg Lie algebra or, more precisely, of its universal enveloping 
algebra. The algebra $L$ is a universal enveloping algebra of a three-dimensional 
\mbox{$\mathbb{Z}_{2}^{3}$-graded} generalized Lie algebra (see Appendix A).
When  anticommutators in the left-hand side of (\ref{rel:colLHeis}) are changed 
into commutators, we indeed have the relations between generators in the universal 
enveloping algebra of the Heisenberg Lie algebra. 

Since the 1970's, generalized (color) Lie algebras have been an object of constant 
interest in both mathematics and physics \cite{Agrav,BahtMPZ,CNS,Gr-Jar},
\cite{Kac-art-adv}--\cite{P}, \cite{R-W(GL)}, \cite{S}--\cite{S-clas}. 
Description of representations of these algebras is an important and interesting 
general problem. It is well known that representations of three-dimensional Lie 
algebras play an important role in the representation theory of general Lie algebras 
and groups, both as test examples and building blocks.
Similarly, one would expect the same to be true for three-dimensional color Lie 
algebras and superalgebras with respect to general color Lie algebras and 
superalgebras. The representations of non-isomorphic  algebras have different
structure. In \cite{S-PhD,S-clas}, three-dimensional color Lie algebras are 
classified in terms of their structure constants, that is in terms of commutation 
relations between generators. In \cite{PSS}, \cite{LS} and \cite{SilSig}, quadratic 
central elements and involutions on these algebras are calculated. In \cite{O-Sil} 
and \cite{S-stm}, Hilbert space $\ast$-representations are described for the graded 
analogues of the Lie algebra $\mathfrak{sl}(2\, ;\,\mathbb{C})$ and of the Lie
algebra of the group of plane motions, two of the non-trivial algebras from the 
classification. The classification of $\ast$-representations in \cite{O-Sil} and 
\cite{S-stm} is achieved, using the method of dynamical systems based on generalized 
Mackey imprimitivity systems.

The graded analogue of the Heisenberg Lie algebra defined by (\ref{rel:colLHeis})
is another important non-trivial algebra in the classification of three-dimensional
color Lie algebras obtained in \cite{S-PhD,S-clas}. In this article we look for 
representations of this algebra. Here, however, we approach representations in a 
totally different way than it was done in \cite{O-Sil} and \cite{S-stm}. 
In this paper we are interested in describing those representations which can be
obtained as power series in representations of Heisenberg's canonical commutation 
relations. 
   
In Section \ref{sec:ser} we show that, with a natural choice for $A_1$ as the first 
generator of the Heisenberg algebra corresponding to differentiation, there are no 
non-zero polynomials in Heisenberg generators which can be taken as $A_2$ and $A_3$ 
so that (\ref{rel:colLHeis}) is satisfied. This means, in particular, that when 
$A_1$ is the differentiation operator, $A_2$ and $A_3$ cannot be chosen as 
differential operators of finite order with polynomial coefficients. We prove 
however that it is possible for $A_2$ and $A_3$ to be power series in the Heisenberg 
generators with infinitely many non-zero terms, thus in particular making possible 
the operator representations by the differential operators of infinite order. In 
Theorem \ref{thm:2rel}, we describe all such formal power series solutions $A_2$ and 
$A_3$ for the first two relations in (\ref{rel:colLHeis}).
In Theorem \ref{thm:3rel}, we present all formal power series solutions $A_2$
and $A_3$ satisfying all three relations in (\ref{rel:colLHeis}).
It turns out that these solutions depend on the choice of two arbitrary odd power 
series, and thus on countably many arbitrary complex parameters. In other words, 
we get two mappings from the sequence space $\mathbb{C}^{\mathbb{N}}$ to formal 
power series in Heisenberg generators, such that elements of their image spaces 
together with the first Heisenberg generator satisfy the commutation relations 
(\ref{rel:colLHeis}). In all these solutions one finds a special series, which 
turns out to be an abstract series generalization of the parity operator, playing
an important role in quantum mechanics, quantum field theory and supersymmetry 
analysis.

By choosing various pairs of operators satisfying the Heisenberg canonical 
commutation relation (\ref{eq:Heisrel}) and substituting them into the power
series obeying (\ref{rel:colLHeis}), one can find large classes of operator 
representations of the commutation relations (\ref{rel:colLHeis}). Section 
\ref{sec:part} is exclusively devoted to examples of such representations. 
Many of these representations, we believe, cannot be reached or classified 
using classical methods based on dynamical systems approach extending Mackey 
imprimitivity systems. We think that these operator representations might have 
significant physical applications. It would be of great interest to investigate 
spectral, structural and analytical properties of such representations on various 
spaces. It also turns out that for some of these representations, the commutation
relations (\ref{rel:colLHeis}) lead to non-trivial functional differential-difference 
interpolation and combinatorial identities involving Euler, Bernoulli and Stirling
numbers.  

\section{Bosonic power series realizations}
\label{sec:ser}
\markboth{Bosonic power series realizations}{Bosonic power series realizations}
\zerocounters

Consider a set $\{A_1,A_2,A_3\}$ in some complex associative algebra with
unit element $I$ satisfying the following commutation relations 
\begin{equation*}
  A_1A_2 + A_2A_1 = A_3, \quad 
  A_1A_3 + A_3A_1 = 0, \quad 
  A_2A_3 + A_3A_2 = 0.
\end{equation*}
It follows immediately that $A_3^2$ commutes with all three elements
$A_1, A_2$ and $A_3$.  
Suppose there exists a non-zero constant $\alpha$ such that $A_3^2=\alpha^2I$.
{}From the first relation we then obtain 
\begin{equation*}
  A_1A_2 A_3+ A_2A_1A_3 = A_3^2 = \alpha^2I,
\end{equation*}
and using that $A_1A_3=-A_3A_1$ by the second relation, we have
\begin{equation*}
  A_1(A_2 A_3) - (A_2A_3)A_1 = \alpha^2I.
\end{equation*}
Let $\hat{A}_2=\alpha^{-1}A_2$ and $\hat{A}_3=\alpha^{-1}A_3$.
Then we can write
\begin{equation}
  A_1(\hat{A}_2\hat{A}_3) - (\hat{A}_2\hat{A}_3)A_1 =I, \label{eq:Heis} 
\end{equation}
showing that $A_1$ and the combination $\hat{A}_2\hat{A}_3$ satisfy the 
Heisenberg canonical commutation relation. By the way this observation implies 
in particular that (\ref{rel:colLHeis}) together with $A_3^2=\alpha^2I,\
\alpha\neq 0$, cannot be satisfied by bounded operators on a Hilbert space or 
even generally by elements in any unital normed algebra, as this is also the 
case for the Heisenberg canonical commutation relation (\ref{eq:Heisrel}) by the 
famous Wintner-Wielandt result \cite{Put,Win-qm,Wie-Qua}.

Assume that we consider $A_1, A_2$ and $A_3$ as elements of the Heisenberg
algebra ${\mathbb C}\langle A,B \rangle/\langle AB-BA-I\rangle$.
Then (\ref{eq:Heis}) suggests that a reasonable Ansatz is to put $A_1=A$ and 
consider the other two generators $A_2$ and $A_3$ as polynomials in $A$ and $B$ 
having coefficients in $\mathbb C$.
Suppose $A_2$ and $A_3$ are any polynomials in $A$ and $B$. Due to the relation 
$AB=I+BA$, it is clear that $A_2$ and $A_3$ can be rewritten as a linear 
combination of monomials with no $B$ to the right of $A$. When a polynomial 
(or a series) in $A$ and $B$ is written in such a way, we say that it is presented 
in its $(B,A)$-normal form. In the Heisenberg algebra, we know that the set of 
ordered monomials $\{B^jA^k \mid j,k\in{\mathbb N}\}$ is linearly independent.
This fact allows one to reduce the problem of equality of two polynomials in $A$ 
and $B$ to checking whether they have the same coefficients when rewritten in 
$(B,A)$-normal form.

We begin with the following theorem, showing that if $A_1=A$, then one is forced 
to work with series in $A$, $B$ with infinitely many non-zero terms, in order to 
be able to find non-trivial realizations of the commutation relations
\begin{equation*}
  A_1A_2 + A_2A_1 = A_3, \quad 
  A_1A_3 + A_3A_1 = 0
\end{equation*}
in terms of the Heisenberg generators $A$ and $B$.   

\begin{theorem}
\label{thm:2rel}
Assume $A$ and $B$ are two elements in an associative algebra over $\mathbb C$ 
with unit element $I$ satisfying the Heisenberg canonical commutation relation 
$AB-BA=I$. 
\begin{itemize}
\item[$(a)$] 
Let $A_1=A$, and both $A_2$ and $A_3$ be polynomials in $A$ and $B$ with 
complex coefficients. Then it follows that the commutation relations 
\begin{equation*}
  A_1A_2 + A_2A_1 = A_3, \quad 
  A_1A_3 + A_3A_1 = 0 
\end{equation*} 
can only be satisfied if $A_2=A_3=0$.

\item[$(b)$]  Let $A_1=A$, and suppose $A_2$ and $A_3$ are formal power series 
in $A$ and $B$ in the $(B,A)$-normal form, i.e.
\begin{equation*}
  A_2=\sum_{j=0}^{\infty}\sum_{k=0}^{\infty} a_{jk}B^jA^k, \quad 
  A_3=\sum_{j=0}^{\infty}\sum_{k=0}^{\infty} \tilde{a}_{jk}B^jA^k,
\end{equation*} 
where the coefficients $a_{jk}, \tilde{a}_{jk} \in {\mathbb C}$.
Then $A_1$, $A_2$ and $A_3$ satisfy the commutation relations 
\begin{equation*}
  A_1A_2 + A_2A_1 = A_3, \quad 
  A_1A_3 + A_3A_1 = 0
\end{equation*} if and only if
\begin{align*}
  A_2 &= T(A,B)V(A)+BT(A,B)W(A), \\
  A_3 &= T(A,B)W(A),
\end{align*} 
where $V(A)$ and $W(A)$ are power series expressions in $A$ with complex 
coefficients, and $T(A,B)$ is given by
\begin{align*}
  T(A,B) &= \sum_{k=0}^{\infty}\frac{(-2)^k}{k!}B^kA^k.
\end{align*}
\end{itemize}
\end{theorem}

\begin{proof}
(a) Any polynomial in $A$ and $B$ can be written in the $(B,A)$-normal form,
and hence, we can assume that 
\begin{equation}
  A_2=\sum_{j=0}^M\sum_{k=0}^Na_{jk}B^jA^k \label{eq:norm1}
\end{equation} 
for some $M,N\in {\mathbb N}$ and $a_{jk}\in {\mathbb C}$.
Eliminating $A_3$ by use of the commutation relations 
and using $A_1=A$, we obtain 
\begin{equation} 
  A_2A^2 + 2AA_2A + A^2A_2 = 0. \label{eq:quad1}
\end{equation} 
Introduce the notation
\begin{equation*} 
  Q(A,B)=A_2A^2 + 2AA_2A + A^2A_2. 
\end{equation*}
The reordering formula
\begin{equation}
  AB^n = B^nA + nB^{n-1} \label{eq:reord1}
\end{equation}
is valid for $n\geq1$ and follows directly from the Heisenberg commutation 
relation (\ref{eq:Heisrel}) by induction on $n$ \cite[p.\:21]{HelSil-book}. \\
By repeated use of (\ref{eq:reord1}), we readily obtain for $n\geq2$ 
\begin{align}
  A^2B^n &= A(AB^n) = A(B^nA + nB^{n-1}) = AB^nA + nAB^{n-1} \notag \\
         &= (B^nA+nB^{n-1})A + n(B^{n-1}A + (n-1)B^{n-2}) \notag \\
         &= B^nA^2+nB^{n-1}A + nB^{n-1}A + (n-1)nB^{n-2} \notag \\
         &= B^nA^2 + 2nB^{n-1}A + (n-1)nB^{n-2}. \label{eq:reord2}
\end{align}
Since the coefficients $a_{jk}$ in (\ref{eq:norm1}) are allowed to be arbitrary 
complex numbers (including zero) one can, without loss of generality, put
\begin{equation*}
  A_2=\sum_{j=0}^N\sum_{k=0}^N a_{jk}B^jA^k, \qquad N\geq2,
\end{equation*} 
and hence 
\begin{equation}
  A_2A^2 = \sum_{j=0}^N\sum_{k=0}^N a_{jk}B^jA^{k+2}. \label{eq:sum1}
\end{equation}
Applying repeatedly the reordering relations (\ref{eq:reord1}) and (\ref{eq:reord2}),
we further obtain
\begin{align}
  AA_2A &= \sum_{j=0}^N\sum_{k=0}^N a_{jk}AB^jA^{k+1} \notag \\
        &= \sum_{k=0}^N a_{0k}A^{k+2} +
           \sum_{j=1}^N\sum_{k=0}^N a_{jk}(B^jA + jB^{j-1})A^{k+1} \notag \\
        &= \sum_{j=0}^N\sum_{k=0}^N a_{jk}B^jA^{k+2} + 
           \sum_{j=1}^N\sum_{k=0}^N ja_{jk}B^{j-1}A^{k+1} \label{eq:sum2}
\end{align}
and
\begin{align}
  A^2A_2 &= \sum_{j=0}^N\sum_{k=0}^N a_{jk}A^2B^jA^k =
            \sum_{k=0}^N a_{0k}A^{k+2} + 
            \sum_{k=0}^N a_{1k}(BA^2 + 2A)A^k \notag \\
         & \qquad + \sum_{j=2}^N\sum_{k=0}^N a_{jk}(B^jA^2 + 2jB^{j-1}A
          + (j-1)jB^{j-2})A^k \notag \\
         &= \sum_{j=0}^N\sum_{k=0}^Na_{jk}B^jA^{k+2} + 
             \sum_{j=1}^N\sum_{k=0}^N2ja_{jk}B^{j-1}A^{k+1} \notag \\
         & \qquad + \sum_{j=2}^N\sum_{k=0}^N(j-1)ja_{jk}B^{j-2}A^k.  \label{eq:sum3}
\end{align}\\
Inserting expressions (\ref{eq:sum1}), (\ref{eq:sum2}) and (\ref{eq:sum3}) into
(\ref{eq:quad1}), we now have
\begin{multline} \raisetag{28pt} \label{eq:reco2}
  Q(A,B) = A_2A^2 + 2AA_2A + A^2A_2 \\
 \begin{aligned}
  &= \sum_{j=0}^N\sum_{k=0}^N4a_{jk}B^jA^{k+2} +
     \sum_{j=1}^N\sum_{k=0}^N4ja_{jk}B^{j-1}A^{k+1} + 
     \sum_{j=2}^N\sum_{k=0}^N(j-1)ja_{jk}B^{j-2}A^k \\
  &= \sum_{j=0}^N\sum_{k=0}^N4a_{jk}B^jA^{k+2} +
     \sum_{j=0}^{N-1}\sum_{k=-1}^{N-1}4(j+1)a_{j+1,k+1}B^jA^{k+2} \\
  &\qquad +
     \sum_{j=0}^{N-2}\sum_{k=-2}^{N-2}(j+1)(j+2)a_{j+2,k+2}B^jA^{k+2} = 0. 
 \end{aligned} 
\end{multline}  
Rearranging the sums, we arrive at
\begin{align*} 
  & Q(A,B)= \sum_{j=0}^{N-2}\sum_{k=0}^{N-2}[4a_{jk} + 4(j+1)a_{j+1,k+1} +
            (j+1)(j+2)a_{j+2,k+2}]B^jA^{k+2} \\
  &\quad \qquad + \sum_{j=0}^{N-2}(j+1)(j+2)a_{j+2,0}B^j +
  \sum_{j=0}^{N-2}(j+1)(j+2)a_{j+2,1}B^jA \\
  &\quad \qquad + \sum_{j=0}^{N-1}4(j+1)a_{j+1,0}B^jA + 
      \sum_{j=0}^{N-1}4(j+1)a_{j+1,N}B^jA^{N+1} \\
  &\quad \qquad + \sum_{k=0}^{N-2}4Na_{N,k+1}B^{N-1}A^{k+2} + 
      \sum_{j=0}^{N}4a_{j,N-1}B^jA^{N+1} +
      \sum_{j=0}^{N}4a_{jN}B^jA^{N+2} \\
  &\quad\qquad + \sum_{k=0}^{N-2}4a_{N-1,k}B^{N-1}A^{k+2} +
      \sum_{k=0}^{N-2}4a_{Nk}B^NA^{k+2}. 
\end{align*}
The polynomial $Q(A,B)$ can now be presented in its $(B,A)$-normal form, 
allowing us to express (\ref{eq:quad1}) as follows
\begin{align}
  & Q(A,B) = \sum_{j=0}^{N-2}\sum_{k=0}^{N-2}[4a_{jk} + 4(j+1)a_{j+1,k+1} +
                   (j+1)(j+2)a_{j+2,k+2}]B^jA^{k+2} \notag \\
  &\qquad \qquad + \sum_{j=0}^{N-2}(j+1)(j+2)a_{j+2,0}B^j  \notag \\
  &\qquad \qquad + \sum_{j=0}^{N-2}[4(j+1)a_{j+1,0} + (j+1)(j+2)a_{j+2,1}]B^jA + 
                   4Na_{N0}B^{N-1}A \notag \\
  &\qquad \qquad + \sum_{j=0}^{N-1}4[a_{j,N-1} + (j+1)a_{j+1,N}]B^jA^{N+1} +  
                   4a_{N,N-1}B^NA^{N+1}  \notag \\
  &\qquad \qquad + \sum_{k=0}^{N-2}4(a_{N-1,k} + Na_{N,k+1})B^{N-1}A^{k+2} + 
                   \sum_{k=0}^{N-2}4a_{Nk}B^NA^{k+2} \notag \\
  &\qquad \qquad + \sum_{j=0}^{N}4a_{jN}B^jA^{N+2} = 0. \label{eq:seco2} 
\end{align}
By linear independence of the set of ordered monomials 
$\{B^jA^k \mid j,k\in{\mathbb N}\}$, all coefficients must be equal to zero, 
giving rise to the following recurrence relation
\begin{equation}
  4a_{jk} + 4(j+1)a_{j+1,k+1} + (j+1)(j+2)a_{j+2,k+2} = 0, \label{eq:rec1}   
\end{equation}
valid for $j,k\in\{0,\ldots,N-2\}$, together with the boundary conditions
\begin{align*}
  (j+1)(j+2)a_{j+2,0} &= 0, \qquad  j=0,\ldots,N-2, \\
  4(j+1)a_{j+1,0} + (j+1)(j+2)a_{j+2,1} &= 0, \qquad
  j=0,\ldots,N-2, \\                      
  Na_{N0} &= 0, \\
  a_{j,N-1} + (j+1)a_{j+1,N} &= 0, \qquad j=0,\ldots,N-1, \\
  a_{N,N-1} &= 0, \\
  a_{N-1,k} + Na_{N,k+1} &= 0,  \qquad k=0,\ldots,N-2, \\
  a_{Nk} &= 0,  \qquad k=0,\ldots,N-2, \\
  a_{jN} &= 0,  \qquad j=0,\ldots,N. 
\end{align*}
It immediately follows, that this system of equations has the solution
\begin{align}
  a_{j0} &= 0, \qquad  j=2,\ldots,N, \notag \\
  a_{21} &= -2a_{10}, \quad a_{j1}= 0, \quad  j=3,\ldots,N, \notag \\
  a_{j,N-1}=a_{jN} &= 0, \qquad j=0,\ldots,N, \label{eq:bousol1} \\            
  a_{N-1,k}= a_{Nk} &= 0, \qquad k=0,\ldots,N. \label{eq:bousol2} 
\end{align}
Consider the square matrix $(a_{ij})$ of size $(N+1)\times(N+1)$.
As expressed by the conditions (\ref{eq:bousol1})--(\ref{eq:bousol2}),
we see that the last two rows and last two columns consist merely
of zeros. In view of relation (\ref{eq:rec1}), it clearly follows 
that $(a_{ij})$ must be the zero matrix, i.e.~all coefficients $a_{jk}=0$, 
showing that (\ref{eq:quad1}) cannot be satisfied by a non-zero polynomial 
expression in the form
$$A_2=\sum_{j=0}^M\sum_{k=0}^Na_{jk}B^jA^k$$
for any $M,N\in {\mathbb N}$, proving part (a) of the theorem.
\\[20pt]
(b) Having now
$$ A_2=\sum_{j=0}^{\infty}\sum_{k=0}^{\infty} a_{jk}B^jA^k,$$
it follows immediately that (\ref{eq:reco2}) will be replaced by 
\begin{multline*} 
  A_2A^2 + 2AA_2A + A^2A_2 \\
 \begin{aligned}
  &= \sum_{j=0}^{\infty}\sum_{k=0}^{\infty}[4a_{jk} + 4(j+1)a_{j+1,k+1} +
            (j+1)(j+2)a_{j+2,k+2}]B^jA^{k+2} \\
  & \quad + \sum_{j=0}^{\infty}(j+1)(j+2)a_{j+2,0}B^j +
     \sum_{j=0}^{\infty}(j+1)[4a_{j+1,0} + (j+2)a_{j+2,1}]B^jA = 0.
 \end{aligned}
\end{multline*}
We still have the recurrence relation 
\begin{equation}
  4a_{jk} + 4(j+1)a_{j+1,k+1} + (j+1)(j+2)a_{j+2,k+2} = 0, \label{eq:rec2}            
\end{equation}
now valid for $j,k\in {\mathbb N}$, together with the conditions
\begin{align*}
  (j+1)(j+2)a_{j+2,0} &= 0, \qquad  j\in {\mathbb N}, \\
  (j+1)[4a_{j+1,0} + (j+2)a_{j+2,1}] &= 0, \qquad j\in {\mathbb N}                  
\end{align*}
or more explicitly
\begin{align}
  a_{j0} &= 0, \qquad  j=2,3,\ldots, \label{eq:bousol3} \\
  a_{21} &= -2a_{10}, \label{eq:bousol4} \\
  a_{j1} &= 0, \qquad  j=3,4,\ldots. \label{eq:bousol5} 
\end{align}
As a consequence of (\ref{eq:rec2}), (\ref{eq:bousol3}) and (\ref{eq:bousol5})
we have
\begin{equation}
  a_{i+2+j,j} = 0, \ \  i,j\in {\mathbb N}. \label{eq:sol1}            
\end{equation}
The elements $a_{j+1,j},\ j\in {\mathbb N}$, can be computed from (\ref{eq:rec2}) 
and (\ref{eq:bousol4}), and will be treated separately below.
In the relation (\ref{eq:rec2}), we now put $k=j+l$, obtaining
\begin{equation}
  4a_{j,j+l} + 4(j+1)a_{j+1,j+l+1} + (j+1)(j+2)a_{j+2,j+l+2} = 0, \label{eq:rec3}      
\end{equation}
where $j,l\in {\mathbb N}$. The substitution $b_j^l=a_{j,j+l}$ brings (\ref{eq:rec3}) 
to the form
\begin{equation}
  4b_j^l + 4(j+1)b_{j+1}^l + (j+1)(j+2)b_{j+2}^l = 0. \label{eq:rec4}
\end{equation}
Suppressing the superscript $l$ for a moment, we look for a general solution 
$(b_{i})_{i=0}^{\infty}$ to the difference equation
\begin{equation}
  4b_j + 4(j+1)b_{j+1} + (j+1)(j+2)b_{j+2} = 0. \label{eq:rec5}
\end{equation}
In order to solve (\ref{eq:rec5}), we introduce the generating function
\begin{equation}
  y(t) = \sum_{j=0}^{\infty}b_jt^j. \label{eq:gen1}
\end{equation}
Successive differentiation of the series expression yields
\begin{align}
  y'(t) &= \sum_{j=0}^{\infty}(j+1)b_{j+1}t^j, \label{eq:gen2} \\
  y''(t) &= \sum_{j=0}^{\infty}(j+1)(j+2)b_{j+2}t^j. \label{eq:gen3} 
\end{align}
(\ref{eq:rec5}), together with the expressions (\ref{eq:gen1}), (\ref{eq:gen2}) 
and (\ref{eq:gen3}), gives rise to the differential equation
\begin{equation}
  y''(t) + 4y'(t) + 4y(t) = 0, \label{eq:gen4}
\end{equation}
where $y$ has to satisfy the initial conditions
\begin{equation}
  y(0) = b_0, \quad y'(0) = b_1. \label{eq:gen5}
\end{equation}
The characteristic equation $(r+2)^2=0$ has the root $r=-2$ of
multiplicity 2, so the general solution has the form
\begin{equation*}
  y(t) = (C_1 +C_2t)e^{-2t}. 
\end{equation*}
By conditions (\ref{eq:gen5}) 
\begin{equation*}
  C_1 = b_0, \quad C_2 = b_1 + 2b_0, 
\end{equation*}
so the solution is
\begin{equation*}
  y(t) = (b_0 + b_1t + 2b_0t)e^{-2t}.
\end{equation*}
Expanding the exponential function yields
\begin{align*} 
  y(t) &= \sum_{j=0}^{\infty}b_jt^j =  
          b_0\sum_{j=0}^{\infty}\frac{(-2)^j}{j!}t^j + 
         (2b_0 + b_1)\sum_{j=0}^{\infty}\frac{(-2)^j}{j!}t^{j+1} \\
       &= b_0 + \sum_{j=0}^{\infty}\left[b_0\frac{(-2)^{j+1}}{(j+1)!} +
          2b_0\frac{(-2)^j}{j!} + b_1\frac{(-2)^j}{j!}\right]t^{j+1}.         
\end{align*}
After identification of the coefficients we have
$$
  b_{j+1} = \frac{(-2)^j}{j!}\left(b_1 + \frac{2j}{j+1}b_0\right), \qquad j = 1,2,\ldots. 
$$
Going back to the earlier notation this means
\begin{equation*}
  b_{j+1}^l = \frac{(-2)^j}{j!}\left(b_1^l + \frac{2j}{j+1}b_0^l\right), \qquad 
  j = 1,2,\ldots,\ \ l = 0,1,\ldots, 
\end{equation*}
and hence,
\begin{equation}
  a_{j+1,j+l+1} = \frac{(-2)^j}{j!}\left(a_{1,l+1} + \frac{2j}{j+1}a_{0l}\right), \quad 
  j = 1,2,\ldots,\ \ l = 0,1,\ldots. \label{eq:gen6}
\end{equation}
If, in the formula (\ref{eq:gen6}) we put $j=1,\ l=-1$, then
\begin{equation*}
  a_{21} = -2(a_{10} + a_{0,-1}) = -2a_{10} 
\end{equation*}
by introducing an auxiliary coefficient $a_{0,-1}=0$.
The solution to the problem (\ref{eq:rec2}), (\ref{eq:bousol3}), (\ref{eq:bousol4}) 
and (\ref{eq:bousol5}) can now be written
\begin{align}
  a_{j+1,j+l} &= \frac{(-2)^j}{j!}\left(a_{1l} 
               + \frac{2j}{j+1}a_{0,l-1}\right), \label{eq:sol2} \\ 
  a_{j+l+2,l} &= 0  \label{eq:sol3},  
\end{align}
where we have $j,l\in {\mathbb N}$.
By virtue of (\ref{eq:sol3}), we may write
\begin{align}
  A_2 &= \sum_{j=0}^{\infty}\sum_{k=0}^{\infty} a_{jk}B^jA^k 
       = \sum_{k=0}^{\infty} a_{0k}A^k  
       + \sum_{j=1}^{\infty}\sum_{k=0}^{\infty} a_{jk}B^jA^k \notag \\
      &= \sum_{k=0}^{\infty} a_{0k}A^k  
         + \sum_{l=0}^{\infty}\sum_{j=0}^{\infty}
         a_{j+1,j+l}B^{j+1}A^{j+l} \label{eq:sol4}
\end{align}
and hence, directly by reordering formula (\ref{eq:reord1})
\begin{align}
  A_3 &= AA_2 + A_2A = \sum_{k=0}^{\infty} a_{0k}A^{k+1}
       + \sum_{l=0}^{\infty}\sum_{j=0}^{\infty}a_{j+1,j+l}AB^{j+1}A^{j+l} \notag \\
      & \qquad + \sum_{k=0}^{\infty} a_{0k}A^{k+1}     
       + \sum_{l=0}^{\infty}\sum_{j=0}^{\infty}a_{j+1,j+l}B^{j+1}A^{j+l+1} \notag \\
      &= 2A_{2}A + \sum_{l=0}^{\infty}\sum_{j=0}^{\infty}(j+1)a_{j+1,j+l}B^jA^{j+l}. 
         \label{eq:sol5}   
\end{align}
Let $c_0 = 0$ and assume  $(c_{i})_{i=1}^{\infty}$ and $(d_{i})_{i=0}^{\infty}$ are 
arbitrary sequences of complex numbers. For $l = 0,1,2,\ldots$ we put 
\begin{equation}
  a_{0,l-1} = c_l, \ \ a_{1l} = d_l \label{eq:cff}
\end{equation}
and have then for $j,l\in {\mathbb N}$
\begin{equation}
  a_{j+1,j+l} = \frac{(-2)^j}{j!}\left(d_l + \frac{2j}{j+1}c_l\right). \label{eq:sol6}
\end{equation}
Inserting now the expressions (\ref{eq:cff}) and (\ref{eq:sol6}) into (\ref{eq:sol4}) 
and (\ref{eq:sol5}), yields
\begin{align*}
  A_2 &= \sum_{k=0}^{\infty}c_{k+1}A^k  
       + \sum_{l=0}^{\infty}\sum_{j=0}^{\infty}
         \frac{(-2)^j}{j!}\left(d_l +
         \frac{2j}{j+1}c_l\right)B^{j+1}A^{j+l} \\
      &= \sum_{k=0}^{\infty}c_{k+1}A^k  
       + \sum_{l=0}^{\infty}\sum_{j=0}^{\infty}
         \frac{(-2)^{j+1}}{(j+1)!}[-(j+1)d_l/2 -
           jc_l]B^{j+1}A^{j+l} \\
      &= \sum_{\substack{k,l=0\\k+l\neq0}}^{\infty}
         \frac{(-2)^k}{k!}[(1-k)c_l - kd_l/2]B^kA^{k+l-1}, 
\end{align*}
\begin{align*}
  A_3 &= 2A_2A + \sum_{l=0}^{\infty}\sum_{j=0}^{\infty}
           \frac{(-2)^j}{j!}(j+1)\left(d_l +
           \frac{2j}{j+1}c_l\right)B^jA^{j+l} \\
      &= 2A_{2}A + \sum_{k=0}^{\infty}\sum_{l=0}^{\infty}
           \frac{(-2)^k}{k!}[(k+1)d_l + 2kc_l]B^kA^{k+l} \\
      &=\sum_{k=0}^{\infty}\sum_{l=0}^{\infty}
           \frac{(-2)^k}{k!}[2(1-k)c_l - kd_l +
           (k+1)d_l + 2kc_l]B^kA^{k+l} \\
      &=\sum_{k=0}^{\infty}\sum_{l=0}^{\infty}
           \frac{(-2)^k}{k!}(2c_l + d_l)B^kA^{k+l}.
\end{align*}
For $l\in {\mathbb N}$, we introduce the coefficients $w_l = 2c_l + d_l$,
allowing us to write
\begin{align*}
  A_2 &= \sum_{\substack{k,l=0\\k+l\neq0}}^{\infty}
         \frac{(-2)^k}{k!}(c_l - \tfrac{1}{2}kw_l)B^kA^{k+l-1}, \\
  A_3 &= \sum_{k=0}^{\infty}\sum_{l=0}^{\infty}
         \frac{(-2)^k}{k!}w_lB^kA^{k+l}, 
\end{align*}
where the coefficients $c_l, w_l$ are arbitrary complex constants except for $c_0$, 
since by definition (\ref{eq:cff}) we have $c_0=0$.
Alternatively, we can express $A_2$ and $A_3$ by separating the summations as follows
\begin{align}
  A_2 &= \sum_{k=0}^{\infty}\sum_{l=1}^{\infty}
          \frac{(-2)^k}{k!}c_lB^kA^{k+l-1} +
          \sum_{k=1}^{\infty}\sum_{l=0}^{\infty}
          \frac{(-2)^{k-1}}{(k-1)!}w_lB^kA^{k-1+l} \notag \\
      &= \sum_{k=0}^{\infty}\frac{(-2)^k}{k!}B^kA^k
           \sum_{l=1}^{\infty}c_lA^{l-1} +
           B\sum_{k=0}^{\infty}\frac{(-2)^k}{k!}B^kA^k
           \sum_{l=0}^{\infty}w_lA^l, \label{eq:A_2sepser} \\
  A_3 &= \sum_{k=0}^{\infty}\frac{(-2)^k}{k!}B^kA^k
           \sum_{l=0}^{\infty}w_lA^l. \label{eq:A_3sepser}
\end{align}
\end{proof}

In the following lemma, we formulate some basic rules, that will be used 
frequently below in proving a corollary to Theorem \ref{thm:2rel} and for the 
proof of our main theorem (Theorem \ref{thm:3rel}).

\begin{lemma}
\label{lem:rel}
Assume $A$ and $B$ are two elements in some complex associative algebra
with unity $I$ satisfying the Heisenberg canonical 
commutation relation $AB-BA=I$. Let $f(A)$ and $g(B)$ be formal power 
series in $A$ and $B$, respectively. Denote by $f'(A)$ and $g'(B)$ their
formal derivatives, obtained by term-wise differentiation of the series
expressions, and let $T(A,B)$ be given by 
\begin{equation}
  T(A,B) = \sum_{k=0}^{\infty}\frac{(-2)^k}{k!}B^kA^k. \label{eq:toper}
\end{equation}
Then the following relations hold true
\begin{itemize}
  \item[$(a)$] $T(A,B)^2=T(A,B)T(A,B) = I,$
  \item[$(b)$] $f(A)B = Bf(A) + f'(A)$, \quad $Ag(B) = g(B)A + g'(B)$,
  \item[$(c)$] $AT(A,B) = -T(A,B)A$, \quad $T(A,B)B = -BT(A,B)$,
  \item[$(d)$] $f(A)T(A,B) = T(A,B)f(-A)$, \quad $T(A,B)g(B) = g(-B)T(A,B)$.
\end{itemize}
\end{lemma}

\begin{proof}
The reordering relation \cite[Cor 2.4 p.\,24]{HelSil-book}
\begin{equation}
  A^iB^j = \sum_{\nu=0}^{\min(i,j)}\nu! 
           \binom{i}{\nu}\binom{j}{\nu}B^{j-\nu}A^{i-\nu} \label{eq:reord3}
\end{equation}
is valid for all non-negative $i$ and $j$, as long as $A$ and $B$ satisfy the 
Heisenberg canonical commutation relation (\ref{eq:Heisrel}). \\
(a) Using the definition (\ref{eq:toper}) and relation (\ref{eq:reord3}), we have
\begin{align*}
  T(A,B)T(A,B) &= \sum_{k=0}^{\infty}\frac{(-2)^k}{k!}B^kA^k
                   \sum_{m=0}^{\infty}\frac{(-2)^m}{m!}B^mA^m \\
               &= \sum_{k=0}^{\infty}\sum_{m=0}^{\infty}
                   \frac{(-2)^{k+m}}{k!\ m!}B^kA^kB^mA^m \\
               &= \sum_{k=0}^{\infty}\sum_{m=0}^{\infty}
                   \sum_{\nu=0}^{\min(k,m)}\frac{(-2)^{k+m}}{k!\ m!}\nu!
                   \binom{k}{\nu}\binom{m}{\nu}B^{k+m-\nu}A^{k+m-\nu} \\
               &= \sum_{k=0}^{\infty}\sum_{m=0}^{\infty}
                   \sum_{\nu=0}^{\min(k,m)}\frac{(-2)^{k+m}}{k!\ (m-\nu)!}
                   \binom{k}{\nu}B^{k+m-\nu}A^{k+m-\nu}.
\end{align*}
Introducing a new summation index $r = k + m -\nu$, this can be expressed as
\begin{equation*}
  T(A,B)T(A,B) = \sum_{r=0}^{\infty}d_rB^rA^r,
\end{equation*} 
where
\begin{align*}
  d_r &= \sum_{k=0}^{r}\sum_{m=r-k}^{r}\frac{(-2)^{k+m}}{k!\ (r-k)!}\binom{k}{k+m-r} \\
      &= (-2)^r\sum_{k=0}^{r}\frac{1}{k!\ (r-k)!}
               \sum_{\nu=0}^{k}(-2)^{\nu}\binom{k}{\nu} \\
      &= \frac{(-2)^r}{r!}\sum_{k=0}^{r}(-1)^k
         \binom{r}{k} = \frac{(-2)^r}{r!}\delta_{r0} = \delta_{r0}.
\end{align*}
Hence,
\begin{equation*}
  T(A,B)T(A,B) = \sum_{r=0}^{\infty}d_rB^rA^r = 
                 \sum_{r=0}^{\infty}\delta_{r0}B^rA^r = I.
\end{equation*} 
(b) From the reordering relation (\ref{eq:reord3}), we obtain  
 $A^nB = BA^n + nA^{n-1}$, valid for all $n\geq1$.
Thus, the first relation in (b) is true for $f(A)=A^n$, $n\in {\mathbb N}$,
and so clearly it also holds for every formal power series in $A$. 
Moreover, for $n\geq1$ we have by (\ref{eq:reord1}) or from relation 
(\ref{eq:reord3}) that $AB^n = B^nA + nB^{n-1}$, showing that the second 
equation (b) holds for $g(B)=B^n$, $n\in {\mathbb N}$,
and therefore it is also valid for any formal power series in $B$. \\
(c) By use of the second relation in (b), we can write
\begin{align*}
  AT(A,B) &= \sum_{k=0}^{\infty}\frac{(-2)^k}{k!}AB^kA^k
           = \sum_{k=0}^{\infty}\frac{(-2)^k}{k!}(B^kA)A^k
           + \sum_{k=1}^{\infty}\frac{(-2)^k}{k!}(kB^{k-1})A^k \\
          &= \sum_{k=0}^{\infty}\frac{(-2)^k}{k!}B^kA^kA
           -2 \sum_{k=1}^{\infty}\frac{(-2)^{k-1}}{(k-1)!}B^{k-1}A^k \\
          & = T(A,B)A - 2T(A,B)A = -T(A,B)A.
\end{align*}
 Applying the first relation in (b), it follows that
\begin{align*}
  T(A,B)B &= \sum_{k=0}^{\infty}\frac{(-2)^k}{k!}B^kA^kB
           = \sum_{k=0}^{\infty}\frac{(-2)^k}{k!}B^k(BA^k)
           + \sum_{k=1}^{\infty}\frac{(-2)^k}{k!}B^k(kA^{k-1}) \\
          &= B\sum_{k=0}^{\infty}\frac{(-2)^k}{k!}B^kA^k
           -2B\sum_{k=1}^{\infty}\frac{(-2)^{k-1}}{(k-1)!}B^{k-1}A^{k-1} \\
          &= BT(A,B) - 2BT(A,B) = -BT(A,B).
\end{align*}
(d) Using the first of relations (c), we have by induction on $n$ that
\begin{equation*}
  A^nT(A,B) = T(A,B)(-A)^n
\end{equation*} 
for all $n\geq0$, and hence, the first relation follows for any power series $f(A)$. \\
{}From the second relation in (c), we have by induction on $n$ that
\begin{equation*}
  T(A,B)B^n = (-B)^nT(A,B)
\end{equation*} 
for all $n\geq0$, and hence, the second relation holds for an arbitrary power series $g(B)$.
\end{proof}

\begin{remark}
\label{rem:pargen}
The series $T(A,B)$ can be seen as an abstract generalization of the parity
operator $f(x) \mapsto f(-x)$. The usual parity operator is obtained in the 
special case of canonical representation of the Heisenberg relation 
(\ref{eq:Heisrel}) when $A=\partial_x:f(x) \mapsto f'(x)$ is differentiation 
and $B=M_x:f(x) \mapsto xf(x)$ is multiplication operator acting on 
functions on $\mathbb R$.
This is proved in the beginning of section \ref{sec:part}.
\end{remark}

In view of the discussion in the beginning of Section \ref{sec:ser},
it is of interest to have a closer look at the compositions $A_2A_3$
and $A_3^2$.

\begin{cor}
\label{cor:2rel}
Assume $A$ and $B$ are two elements in some complex associative algebra 
with unity $I$ satisfying the Heisenberg canonical 
commutation relation $AB-BA=I$. Let $A_1=A$, and suppose $A_2$ and $A_3$ are
formal power series in $A$ and $B$ in the $(B,A)$-normal form, i.e.
\begin{equation*}
  A_2=\sum_{j=0}^{\infty}\sum_{k=0}^{\infty} a_{jk}B^jA^k, \quad
  A_3=\sum_{j=0}^{\infty}\sum_{k=0}^{\infty} \tilde{a}_{jk}B^jA^k,
\end{equation*} 
where the coefficients $a_{jk}, \tilde{a}_{jk} \in {\mathbb C}$.
If $A_1$, $A_2$ and $A_3$ satisfy the commutation relations 
\begin{equation*}
  A_1A_2 + A_2A_1 = A_3, \quad 
  A_1A_3 + A_3A_1 = 0
\end{equation*} 
then it follows that
\begin{align*}
  A_2A_3 & = V(-A)W(A)+BW(-A)W(A), \\
  A_3^2 & = W(-A)W(A),
\end{align*} 
where $V(A)$ and $W(A)$ are formal power series expressions in $A$
with complex coefficients.
\end{cor}

\begin{proof}
The general solution to the problem according to Theorem \ref{thm:2rel} 
can be expressed as 
\begin{align}
  A_2 &= T(A,B)V(A)+BT(A,B)W(A), \label{eq:exprA_2} \\
  A_3 &= T(A,B)W(A), \label{eq:exprA_3}
\end{align}
where $V(A)$ and $W(A)$ are arbitrary formal power series with
coefficients in $\mathbb C$.
Direct substitution of these expressions and application of the rules
in Lemma \ref{lem:rel}, yields
\begin{align*}
  A_2A_3 &= T(A,B)V(A)T(A,B)W(A)+BT(A,B)W(A)T(A,B)W(A)   \\
             &= T(A,B)T(A,B)V(-A)W(A)+BT(A,B)T(A,B)W(-A)W(A) \\
             &= V(-A)W(A)+BW(-A)W(A), \\
  A_3A_3 &= T(A,B)W(A)T(A,B)W(A) = T(A,B)T(A,B)W(-A)W(A) \\
             & = W(-A)W(A). 
\end{align*}
\end{proof}

In Theorem \ref{thm:2rel} we stated the general $(B,A)$-normal form of 
the power series $A_2$ and $A_3$ satisfying the first two relations in 
(\ref{rel:colLHeis}). We now turn our attention to investigating the 
possibility of satisfying also the third relation, namely
\begin{equation}
  A_2A_3 + A_3A_2 = 0. \label{eq:rel3}
\end{equation}

In the following theorem, being the main result of this article, we give the 
general solution to the problem with all three relations (\ref{rel:colLHeis}). 
In the formulation the exponential generating function $E(t)$ of the so-called
Euler numbers is used. For basic facts about the Euler polynomials and Euler 
numbers, we refer to Appendix B.

\begin{theorem}
\label{thm:3rel}
Suppose $A$ and $B$ are two elements in a unital associative 
algebra over $\mathbb C$ with unity $I$ satisfying the Heisenberg canonical 
commutation relation $AB-BA=I$. Put $A_1=A$ and let $A_2$, $A_3$ and $T(A,B)$ be
formal power series in $A$ and $B$ in the $(B,A)$-normal form given as
\begin{equation*}
  A_2=\sum_{j,k=0}^{\infty} a_{jk}B^jA^k, \ \ 
  A_3=\sum_{j,k=0}^{\infty}\tilde{a}_{jk}B^jA^k, \ \ 
  T(A,B) = \sum_{k=0}^{\infty}\frac{(-2)^k}{k!}B^kA^k,
\end{equation*} 
where the coefficients $a_{jk},\tilde{a}_{jk}\in {\mathbb C}$.
If $A_1$, $A_2$ and $A_3$ satisfy the commutation relations 
\begin{equation*}
  A_1A_2 + A_2A_1 = A_3,  \quad 
  A_1A_3 + A_3A_1 = 0,  \quad
  A_2A_3 + A_3A_2 = 0,
\end{equation*}
then either $A_3 = 0$ and $A_2 = T(A,B)V(A)$,
where $V(A)$ is a formal power series in $A$ with complex coefficients, or
\begin{align*}
  A_2 &=c\,T(A,B)E(\varphi(A))[e^{\varphi(A)}\psi(A)- \tfrac{1}{2}\varphi'(A)]
        +c\,BT(A,B)e^{\varphi(A)},\\
  A_3 &= c\,T(A,B)e^{\varphi(A)},
\end{align*} 
where $c$ is a non-zero complex constant, $E(t)$ is the exponential
generating function of the Euler numbers $E_k$ given by
\begin{equation*}
  E(t)= \sum_{n=0}^{\infty}\frac{E_{2n}}{(2n)!}t^{2n},
\end{equation*}
and both $\varphi(A)$ and $\psi(A)$ are odd formal power series expressions in $A$ 
with complex coefficients.
\end{theorem}

\begin{proof}
By Theorem \ref{thm:2rel} we have, when considering only the first two relations, 
a general solution of the form given by (\ref{eq:exprA_2}) and (\ref{eq:exprA_3}).
In the present case, $A_2$ and $A_3$ are supposed to satisfy the additional 
condition $A_2A_3+A_3A_2=0$. By Corollary \ref{cor:2rel} we have
\begin{equation*}
  A_2A_3 = V(-A)W(A)+BW(-A)W(A). 
\end{equation*}
Inserting the expressions (\ref{eq:exprA_3}) and (\ref{eq:exprA_2}) and applying 
the rules in Lemma \ref{lem:rel}, yields
\begin{align*}
  A_3A_2 &= T(A,B)W(A)T(A,B)V(A)+T(A,B)W(A)BT(A,B)W(A) \\
         &= T(A,B)T(A,B)W(-A)V(A)-T(A,B)W(A)T(A,B)BW(A) \\
         &= W(-A)V(A)-T(A,B)T(A,B)W(-A)BW(A) \\
         &= W(-A)V(A)-[BW(-A)-W'(-A)]W(A) \\
         &= W(-A)V(A)-BW(-A)W(A)+W'(-A)W(A), 
\end{align*}
and hence, we obtain
\begin{equation*}
  A_2A_3+ A_3A_2 = V(-A)W(A)+ W(-A)V(A)+W'(-A)W(A)=0. 
\end{equation*}
We shall now consider the functional-differential equation
\begin{equation}
  V(-A)W(A)+ W(-A)V(A)+W'(-A)W(A)=0.  \label{ode:VW+} 
\end{equation}
{}From $AB-BA=I$ it follows that $(-A)(-B)-(-B)(-A)=I$, and hence together
with $\tilde{A}_1=-A$ the expressions
\begin{align*}
  \tilde{A}_2 &= T(-A,-B)V(-A)-BT(-A,-B)W(-A), \\
  \tilde{A}_3 &= T(-A,-B)W(-A), 
\end{align*}
will satisfy the same commutation relations as (\ref{eq:exprA_2}) and
(\ref{eq:exprA_3}). Thus, it follows that
\begin{equation*}
  \tilde{A}_2\tilde{A}_3+\tilde{A}_3\tilde{A}_2 = V(A)W(-A)+ W(A)V(-A)+W'(A)W(-A)=0,
\end{equation*}
and we also have the equation
\begin{equation}
  V(A)W(-A)+ W(A)V(-A)+W'(A)W(-A)=0,  \label{ode:VW-} 
\end{equation}
We shall look for the general solution to (\ref{ode:VW+}) and
(\ref{ode:VW-} ) in the form of formal power series
\begin{equation*}
  V(A)=\sum_{l=0}^{\infty}v_lA^l, \qquad W(A)=\sum_{l=0}^{\infty}w_lA^l,
       \qquad v_l,w_l\in {\mathbb C}.
\end{equation*}
Subtracting (\ref{ode:VW+}) from (\ref{ode:VW-}) yields the equation
\begin{equation}
  W'(A)W(-A)-W'(-A)W(A)=0. \label{ode:W} 
\end{equation}
Integrating (\ref{ode:W}) and noting that $W(0)=w_0$, we have
\begin{equation}
  W(A)W(-A)=w_0^2. \label{feq:W} 
\end{equation}
Considering first the case when $w_0=0$, one has, by substitution of the 
series $W(A)$ into (\ref{feq:W}), the infinite system of equations
\begin{equation*}  
  \sum_{l=0}^s(-1)^lw_lw_{s-l} = 0, \qquad  s=0,1,2,\ldots,     
\end{equation*}
having the unique solution $w_0=w_1=w_2=\ldots=0$, so that $W(A)=0$.
Putting $W(A)=0$ in the differential equation (\ref{ode:VW+}), we see that
there is no equation left for $V(A)$, i.e.~$V(A)$ can be chosen arbitrarily
in the expression (\ref{eq:exprA_2}) for $A_2$,
proving the first statement of the theorem.

Assuming that $w_0\neq0$, we can divide both sides of equation (\ref{feq:W})
by the non-zero constant $w_0^2$, obtaining the simple equation $g(A)g(-A)=1$,
where $g(A):=W(A)/w_0$. Taking the logarithm of both sides, it follows that 
$\log g(A)$ has to be an odd power series expression since $\log g(A)+\log g(-A)=0$. 
Let $\varphi(A)=\log g(A)$ and we have $W(A)=w_0g(t)=w_0\exp(\varphi(A))$, which 
is the general solution to (\ref{feq:W}), $\varphi(A)$ being any odd formal power 
series with complex coefficients.
Substituting for $W(A)$ into equation (\ref{ode:VW-}) yields
\begin{equation*}
  V(A)w_0\exp(\varphi(-A))+w_0\exp(\varphi(A))V(-A)+w_0^2\varphi'(A)=0.
\end{equation*}
This can be written as
\begin{equation*}
  V(A)[\cosh\varphi(A)-\sinh\varphi(A)]+V(-A)[\cosh\varphi(A)+\sinh\varphi(A)]+w_0\varphi'(A)=0,
\end{equation*}
\begin{equation*}
  [V(A)+V(-A)]\cosh\varphi(A)=[V(A)-V(-A)]\sinh\varphi(A)-w_0\varphi'(A),
\end{equation*}
and we have
\begin{equation*}
  2V(A)\cosh\varphi(A)=[V(A)-V(-A)]\exp(\varphi(A))-w_0\varphi'(A).
\end{equation*}
Denoting the odd part of $V(A)$ by $V_1(A)$, we obtain
\begin{equation}
  V(A)\cosh\varphi(A)=V_1(A)\exp(\varphi(A))-\tfrac{w_0}{2}\varphi'(A). \label{eq:Vcosh}
\end{equation}
Let $E(X)$ denote the formal power series defined by
\begin{equation*}
  E(X):= \sum_{n=0}^{\infty}\frac{E_{2n}}{(2n)!}X^{2n},
         \quad E_k=\rm{Euler\ numbers},
\end{equation*}
being the inverse of the formal power series given by $\cosh(X)$, in the
sense that $E(X)\cosh(X)=\cosh(X)E(X)=1$.
Multiplying both sides of equation (\ref{eq:Vcosh}) by $E(\varphi(A))$ yields
\begin{equation*}
  V(A)=E(\varphi(A))V_1(A)\exp(\varphi(A))-\tfrac{w_0}{2}E(\varphi(A))\,\varphi'(A). 
\end{equation*}
We have an expression for $V(A)$ in terms of the odd power series $\varphi(A)$ and the odd 
part $V_1(A)$ of $V(A)$. $V_1(A)$ can be chosen arbitrarily from the set of formal odd power
series with coefficients from ${\mathbb C}$. Writing this as $V_1(A)=w_0\psi(A)$ with $w_0$
and $\psi(A)$ arbitrary, we have
\begin{equation}
  V(A)=w_0E(\varphi(A))\exp(\varphi(A))\,\psi(A)-\tfrac{w_0}{2}E(\varphi(A))\,\varphi'(A). \label{eq:V}
\end{equation}
where $\varphi(A)$ and $\psi(A)$ are arbitrary odd formal power series with 
complex coefficients.
\end{proof}

The simple expressions obtained for the combinations $A_2A_3$ and $A_3^2$ in 
Corollary \ref{cor:2rel}, being essentially pure series in $A$, can now be given in 
a more explicit form. We formulate the following:
\begin{cor}
\label{cor:3rel}
Suppose $A$ and $B$ are two elements in some complex associative algebra 
with unity $I$ satisfying the Heisenberg canonical commutation relation 
$AB-BA=I$. Put $A_1=A$ and let $A_2$ and $A_3$ be formal power series in 
$A$ and $B$ in the $(B,A)$-normal form given as
\begin{equation*}
  A_2=\sum_{j=0}^{\infty}\sum_{k=0}^{\infty} a_{jk}B^jA^k, \ \ 
  A_3=\sum_{j=0}^{\infty}\sum_{k=0}^{\infty}\tilde{a}_{jk}B^jA^k, 
\end{equation*} 
where the coefficients $a_{jk},\tilde{a}_{jk}\in {\mathbb C}$.
If $A_1$, $A_2$ and $A_3$ satisfy the commutation relations 
\begin{equation*}
  A_1A_2 + A_2A_1 = A_3, \quad 
  A_1A_3 + A_3A_1 = 0, \quad
  A_2A_3 + A_3A_2 = 0,
\end{equation*}
then we have
\begin{align*}
  A_2A_3 & = c\,B-c\,E(\varphi(A))[\psi(A)+\tfrac{1}{2}e^{\varphi(A)}\varphi'(A)], \\
  A_3^2 & = c\,I,
\end{align*}
where $c\in{\mathbb C}$, $E(t)$ is the exponential generating function
of the Euler numbers $E_k$ given by
\begin{equation*}
  E(t)= \sum_{n=0}^{\infty}\frac{E_{2n}}{(2n)!}t^{2n},
\end{equation*}
and both $\varphi(A)$ and $\psi(A)$ are odd formal power series expressions in 
$A$ with complex coefficients.
\end{cor}

\begin{proof}
By Corollary \ref{cor:2rel} we know that
\begin{align*}
  A_2A_3 &= V(-A)W(A)+BW(-A)W(A),\\
  A_3^2 &= W(-A)W(A),
\end{align*} 
where in the present case the formal power series $V(A)$ and $W(A)$,
appearing in the expressions for $A_2$ and $A_3$ in Theorem \ref{thm:2rel}, 
must be chosen such that the third relation $A_2A_3+A_3A_2=0$ is satisfied. 
The general solution to that problem is found in the proof of Theorem
\ref{thm:3rel} and is of the form
\begin{align*}
  V(A) &= w_0E(\varphi(A))\exp(\varphi(A))\,\psi(A)-\tfrac{1}{2}w_0E(\varphi(A))\,\varphi'(A), \\
  W(A) &= w_0\exp(\varphi(A)),
\end{align*}
where  $w_0\in{\mathbb C}$, $\varphi(A)$ and $\psi(A)$ are arbitrary odd
formal power series with coefficients in ${\mathbb C}$ and 
$E(t)=\sum_{n=0}^{\infty}\frac{E_{2n}}{(2n)!}t^{2n}$. This yields
\begin{align*}
  W(-A)W(A) &= w_0^2\exp(\varphi(-A))\exp(\varphi(A))=w_0^2\exp(\varphi(A)+\varphi(-A)) \\
            &= w_0^2\exp(\varphi(A)-\varphi(A))=w_0^2\,I, \\
  V(-A)W(A) &= w_0^2E(\varphi(-A))\exp(\varphi(-A))\,\psi(-A)\exp(\varphi(A)) \\
            & \qquad -\tfrac{1}{2}w_0^2E(\varphi(-A))\,\varphi'(-A)\exp(\varphi(A)) \\
            &= -w_0^2E(\varphi(A))\psi(A)-\tfrac{1}{2}w_0^2E(\varphi(A))\,\varphi'(A)\exp(\varphi(A)) \\
            &= -w_0^2E(\varphi(A))[\psi(A)+\tfrac{1}{2}\exp(\varphi(A))\,\varphi'(A)],
\end{align*}
proving the statement of the corollary.
\end{proof}

\begin{remark}
 Note that exchange of $A_1$ and  $A_2$ does not change commutation relations 
(\ref{rel:colLHeis}). So, by exchanging $A_1$ and $A_2$
in all statements of the article, we get other expressions 
for $A_1$ and $A_2$ in terms of the Heisenberg generators.
\end{remark}

\section{Some particular bosonic representations}
\label{sec:part}
\markboth{Some particular bosonic representations}
{Some particular bosonic representations}
\zerocounters
In this section we will describe some non-trivial particular 
representations of the color analogue of the Heisenberg Lie algebra 
defined by the commutation relations (\ref{rel:colLHeis}). 
All examples are based on the general statement in Theorem \ref{thm:3rel} 
and correspond to simple specific choices of the odd formal power series 
$\varphi(A)$ and $\psi(A)$.
The constant $c$ is unimportant (except when $c=0$, see
Example \ref{exa:A_3=0}) and we will put c=1 unless stated otherwise.

As a concrete example of elements satisfying the Heisenberg commutation 
relation (\ref{eq:Heisrel}), we can consider the differentiation and multiplication
operators $\partial_x$ and $M_x$ defined on the linear space ${\mathbb C}[x]$,
consisting of all complex-valued polynomial functions of a single real variable
$x$. If $f(x)=\sum_{k=0}^nf_kx^k$, then by definition 
 $$(\partial_xf)(x)=\sum_{k=1}^nf_kkx^{k-1},\qquad (M_xf)(x)=xf(x)$$
and we have the well-known relation $\partial_xM_x-M_x\partial_x=I$.
As a basis for ${\mathbb C}[x]$ we can take the set
of monomials $\{1,x,x^2,\ldots,x^n,x^{n+1},\ldots\}$.
Acting on an arbitrary basis vector $x^n$, we find 
\begin{align*}
  T(\partial_x,M_x)(x^n) &= \sum_{k=0}^{\infty}
            \frac{(-2)^k}{k!}M_x^k\partial_x^{k}x^n =
            \sum_{k=0}^n\frac{(-2)^k}{k!}x^k\frac{n!}{(n-k)!}x^{n-k}\\
           &=x^n\sum_{k=0}^n(-2)^k\binom{n}{k}=(-1)^nx^n=(-x)^n.
\end{align*}
In fact, if $g$ is an analytic function on $\mathbb R$, we have 
by Taylor's Theorem
\begin{align}
  e^{\partial_x}g(x) &= \sum_{n=0}^{\infty}\frac{\partial_x^n}{n!}g(x)
                    =\sum_{n=0}^{\infty}\frac{g^{(n)}(x)}{n!}=g(x+1), \label{eq:expder} \\
  T(\partial_x,M_x)g(x)  
  &=\sum_{k=0}^{\infty}\frac{(-2)^k}{k!}M_x^k\partial_x^{k}g(x)
   =\sum_{k=0}^{\infty}\frac{(-2)^k}{k!}x^kg^{(k)}(x) \notag \\
   & =\sum_{k=0}^{\infty}\frac{g^{(k)}(x)}{k!}(-2x)^k=g(x-2x)=g(-x). \label{eq:parity}
\end{align}

\begin{example}
\label{exa:A_3=0}
By Theorem \ref{thm:3rel} we have a solution corresponding to $A_3=0$ given by 
\begin{equation*}
  A_1=A,\quad  A_2=T(A,B)V(A),\quad A_3=0,
\end{equation*}
where $V(A)$
can be any power series in $A$ having complex coefficients.
In this case, there is only one non-trivial relation to 
satisfy. Applying rule (c) of Lemma \ref{lem:rel} we readily obtain
\begin{align*}
  A_1A_2+A_2A_1 &= AT(A,B)V(A)+T(A,B)V(A)A \\
                &= -T(A,B)AV(A)+T(A,B)V(A)A=0.
\end{align*}
For any non-zero $V(A)$, we obviously have $A_2$
given as an infinite power series expression in $A$ and $B$.
This is to be expected as a consequence of Theorem \ref{thm:2rel} (a).
Taking $V(A)=0$, we obtain the trivial realization
$A_1=A$, and $A_2=A_3=0$. As we have shown in Section \ref{sec:ser},
this is the only possible solution, when $A_2$ and $A_3$ 
are polynomials in $A$ and $B$.

Assuming that
$V(A)=\sum_{l=0}^{\infty}v_lA^l$, we obtain for our simple solution above
\begin{align*}
  A_1 &= \partial_x, \\
  A_2 &= T(\partial_x,M_x)V(\partial_x) = \sum_{k=0}^{\infty}
               \frac{(-2)^k}{k!}M_x^k\partial_x^k\sum_{l=0}^{\infty}v_l\partial_x^l, \\
  A_3 &= 0.
\end{align*}
Acting on an arbitrary basis vector $x^n$, we find 
\begin{align*}
  A_1(x^n) &= \partial_xx^n = nx^{n-1}, \\
  V(\partial_x)(x^n) &=\sum_{l=0}^{\infty}v_l\partial_x^lx^n= 
                       \sum_{l=0}^nv_l\frac{n!}{(n-l)!}x^{n-l}, \\
  A_2(x^n) &= T(\partial_x,M_x)V(\partial_x)x^n 
            = \sum_{l=0}^nv_l\,l!\binom{n}{l}(-x)^{n-l}.
\end{align*}
\end{example}

\begin{example}
\label{exa:A_3=T}
Choosing $\varphi=\psi\equiv 0$ in the general solution expressed in 
Theorem \ref{thm:3rel}, we obtain 
\begin{equation*}
  A_1=A,\quad A_2=BT(A,B),\quad A_3=T(A,B)
\end{equation*}
so in this case we have the simple relation $A_2=BA_3$.
Considering the same situation as in Example \ref{exa:A_3=0} with 
$A=\partial_x$ and $B=M_x$ defined on the linear space ${\mathbb C}[x]$, 
we have here
\begin{align*}
  A_1 &= \partial_x, \\
  A_2 &= M_xT(\partial_x,M_x) = \sum_{k=0}^{\infty}
         \frac{(-2)^k}{k!}M_x^{k+1}\partial_x^k, \\
  A_3 &= T(\partial_x,M_x) = \sum_{k=0}^{\infty}
         \frac{(-2)^k}{k!}M_x^k\partial_x^{k}. 
\end{align*}
These operators are defined on the whole of ${\mathbb C}[x]$, and by 
Theorem \ref{thm:3rel} they satisfy relations (\ref{rel:colLHeis}).
Acting on an arbitrary basis vector $x^n$, we obtain
\begin{align*}
  A_1(x^n) &= \partial_xx^n = nx^{n-1}, \\
  A_2(x^n) &= M_x T(\partial_x,M_x)x^n= M_x(-x)^n=(-1)^nx^{n+1}, \\
  A_3(x^n) &= T(\partial_x,M_x)x^n = (-x)^n.
\end{align*}
So, for any polynomial $p(x)\in {\mathbb C}[x]$, we have 
\begin{equation*}
  (A_1p)(x) = p'(x), \quad
  (A_2p)(x) = xp(-x), \quad
  (A_3p)(x) = p(-x).
\end{equation*}
These three operators can be defined for any differentiable function 
$f$ and they satisfy the commutation relations (\ref{rel:colLHeis}), 
as proved by the following simple calculations:
\begin{align*}
  A_1A_2f(x) &= \partial_x(xf(-x))= f(-x)-xf'(-x), \\
  A_2A_1f(x) &= A_2f'(x)=xf'(-x), \\
  A_1A_3f(x) &= \partial_xf(-x)=-f'(-x), \\
  A_3A_1f(x) &= A_3f'(x)=f'(-x), \\
  A_2A_3f(x) &= A_2f(-x)=xf(x), \\
  A_3A_2f(x) &= A_3xf(-x)=-xf(x). 
\end{align*}
So,
\begin{align*}
  (A_1A_2+A_2A_1)f(x) &= f(-x) = A_3f(x), \\
  (A_1A_3+A_3A_1)f(x) &= -f'(-x)+f'(-x)=0, \\
  (A_2A_3+A_3A_2)f(x) &= xf(x)-xf(x)=0. 
\end{align*}
\end{example}

\begin{example}
\label{exa:dirsum}
Let $C^{\infty}({\mathbb R})$ be the set of all complex-valued infinitely 
differentiable functions on the real line. An arbitrary function $f$ can 
be written in a unique way as a sum $f(x)=f_0(x)+ f_1(x)$ of its even and 
odd parts, where
\begin{equation*}
  f_0(x)=\frac{f(x)+ f(-x)}{2}\,, \qquad f_1(x)=\frac{f(x)- f(-x)}{2}\,.
\end{equation*}
If $C^{\infty}_0({\mathbb R})$ and  $C^{\infty}_1({\mathbb R})$ denote the 
subsets of $C^{\infty}({\mathbb R})$ consisting of even and odd infinitely
differentiable functions respectively, then this means that 
$C^{\infty}({\mathbb R})$ can be expressed as a direct sum 
$C^{\infty}({\mathbb R}) = C^{\infty}_0({\mathbb R}) \oplus C^{\infty}_1({\mathbb R})$.
Now, let $A_1$, $A_2$ and $A_3$ be defined on $C^{\infty}({\mathbb R})$, as in the 
previous example, by the equations
\begin{equation*}
  (A_1f)(x) = f'(x), \quad
  (A_2f)(x) = xf(-x), \quad
  (A_3f)(x) = f(-x).
\end{equation*}
By considering different domains of definition for these operators by
restricting to the subspaces considered above, we have
\begin{align*}
  A_1:C^{\infty}({\mathbb R}) \rightarrow C^{\infty}({\mathbb R}), & \quad
  A_1:C^{\infty}_0({\mathbb R}) \rightarrow C^{\infty}_1({\mathbb R}), & 
  A_1:C^{\infty}_1({\mathbb R}) \rightarrow C^{\infty}_0({\mathbb R}), & \\
  A_2:C^{\infty}({\mathbb R}) \rightarrow C^{\infty}({\mathbb R}), & \quad
  A_2:C^{\infty}_0({\mathbb R}) \rightarrow C^{\infty}_1({\mathbb R}), & 
  A_2:C^{\infty}_1({\mathbb R}) \rightarrow C^{\infty}_0({\mathbb R}), & \\
  A_3:C^{\infty}({\mathbb R}) \rightarrow C^{\infty}({\mathbb R}), & \quad
  A_3:C^{\infty}_0({\mathbb R}) \rightarrow C^{\infty}_0({\mathbb R}), &
  A_3:C^{\infty}_1({\mathbb R}) \rightarrow C^{\infty}_1({\mathbb R}).
\end{align*}
where $A_2=M_x$ on $C^{\infty}_0({\mathbb R})$, $A_2=-M_x$ on 
$C^{\infty}_1({\mathbb R})$, $A_3=I$ on $C^{\infty}_0({\mathbb R})$, and
$A_3=-I$ on $C^{\infty}_1({\mathbb R})$.
Now, define $$ \mathbf{A}_i:C^{\infty}_0({\mathbb R}) \oplus C^{\infty}_1({\mathbb R})
  \rightarrow  C^{\infty}_0({\mathbb R}) \oplus C^{\infty}_1({\mathbb R}),\ i=1,2,3$$
where the operators $\mathbf{A}_1$, $\mathbf{A}_2$ and $\mathbf{A}_3$ are given
by the operator matrices
\begin{equation*}
  \mathbf{A}_1=\begin{pmatrix}0 & \partial_x \\ \partial_x & 0 \end{pmatrix}, \quad
  \mathbf{A}_2=\begin{pmatrix}0 & -M_x \\ M_x & 0 \end{pmatrix}, \quad
  \mathbf{A}_3=\begin{pmatrix}I & 0 \\ 0 & -I \end{pmatrix}.
\end{equation*}
Using the relation $\partial_xM_x-M_x\partial_x=I$, one easily verifies by direct
matrix multiplication that the operators $\mathbf{A}_1$, $\mathbf{A}_2$ and 
$\mathbf{A}_3$ satisfy the relations 
\begin{equation*}
  \mathbf{A}_1\mathbf{A}_2+\mathbf{A}_2\mathbf{A}_1 = \mathbf{A}_3, \quad 
  \mathbf{A}_1\mathbf{A}_3+\mathbf{A}_3\mathbf{A}_1 = \mathbf{0},  \quad 
  \mathbf{A}_2\mathbf{A}_3+\mathbf{A}_3\mathbf{A}_2 = \mathbf{0},
\end{equation*} 
where $\mathbf{0}$ denotes the $2\times2$ zero matrix.

Actually, if we let ${\mathbb C}^{\mathbb R}$ and $D({\mathbb R})$ denote the sets 
of all complex-valued functions on ${\mathbb R}$ and all complex-valued 
differentiable functions on ${\mathbb R}$ respectively, then clearly $\mathbf{A}_1$ 
can be defined on the space $D({\mathbb R})=D_0({\mathbb R}) \oplus D_1({\mathbb R})$ 
while $\mathbf{A}_2$ and $\mathbf{A}_3$ are well-defined on the whole of 
${\mathbb C}^{\mathbb R}={\mathbb C}^{\mathbb R}_0 \oplus {\mathbb C}^{\mathbb R}_1$. 
Here subscripts $0$ and $1$ have the same meaning as above, i.e.~indicate subsets of 
even and odd functions. It follows that the three relations will be satisfied if the
domain of definition is chosen as 
$D({\mathbb R})=D_0({\mathbb R}) \oplus D_1({\mathbb R})$. 
\end{example}

\begin{example}
\label{exa:tensor}
The representations of (\ref{rel:colLHeis}) described in Example \ref{exa:dirsum}
can be generalized as follows.
Let $H$ be a linear space and $H_0$ and $H_1$ be subspaces of $H$ such that 
$H_0 \cap H_1 = \{0\}$. Consider $H_0 \oplus H_1$, the subspace of H
which is a direct sum of $H_0$ and $H_1$. Any linear operator $Y$ on
$H_0 \oplus H_1$ can be written as an operator matrix $Y=\left(
\begin{smallmatrix} Y_{00}&Y_{01}\\Y_{10}&Y_{11}\end{smallmatrix} \right)$,
where $Y_{jk}:H_k \rightarrow H_j$ for $j,k \in \{0,1\}$ are linear operators.
Suppose $A$ and $B$ are linear operators on $H$ satisfying on $H$ the Heisenberg
canonical commutation relation $AB-BA=I$. Then the linear operators
\begin{equation*}
  A_1=\begin{pmatrix}0 & A \\ A & 0 \end{pmatrix}, \quad
  A_2=\begin{pmatrix}0 & -B \\ B & 0 \end{pmatrix}, \quad
  A_3=\begin{pmatrix}I & 0 \\ 0 & -I \end{pmatrix}
\end{equation*}
on  $H_0 \oplus H_1$ satisfy the commutation relations (\ref{rel:colLHeis}) of 
the graded analogue of the Heisenberg Lie algebra. 
This can be proved by the following calculations:
\begin{align*}
  A_1A_2+A_2A_1 &= \begin{pmatrix}AB & 0 \\ 0 & -AB \end{pmatrix}
                 + \begin{pmatrix}-BA & 0 \\ 0 & BA \end{pmatrix} 
                 = \begin{pmatrix}I & 0 \\ 0 & -I \end{pmatrix} = A_3,\\
  A_1A_3+A_3A_1 &= \begin{pmatrix}0 & -A \\ A & 0 \end{pmatrix}
                 + \begin{pmatrix}0 & A \\ -A & 0 \end{pmatrix}
                 = \begin{pmatrix}0 & 0 \\ 0 & 0 \end{pmatrix} =  0, \\
  A_2A_3+A_3A_2 &= \begin{pmatrix}0 & B \\ B & 0 \end{pmatrix}
                 + \begin{pmatrix}0 & -B \\ -B & 0 \end{pmatrix}
                 = \begin{pmatrix}0 & 0 \\ 0 & 0 \end{pmatrix} =  0. 
\end{align*}
Another way to form such block representations is to use the tensor product.
For some linear space $H$ we consider the tensor product ${\mathbb C}^2 \otimes H$,
being also a linear space over ${\mathbb C}$. Let $A$ and $B$ be operators on H
satisfying $AB-BA=I$. Then the operators
\begin{equation*}
  A_1=\begin{pmatrix}0 & 1 \\ 1 & 0 \end{pmatrix} \otimes A, \quad
  A_2=\begin{pmatrix}0 & -1 \\ 1 & 0 \end{pmatrix} \otimes B, \quad
  A_3=\begin{pmatrix}1 & 0 \\ 0 & -1 \end{pmatrix} \otimes I
\end{equation*}
on ${\mathbb C}^2 \otimes H$ satisfy relations (\ref{rel:colLHeis}).
Direct computation, using the rules for the tensor product, shows that 
\begin{align*}
  A_1A_2+A_2A_1 &= \begin{pmatrix}0 & 1 \\ 1 & 0 \end{pmatrix}
                   \begin{pmatrix}0 & -1 \\ 1 & 0 \end{pmatrix}\otimes AB 
                 + \begin{pmatrix}0 & -1 \\ 1 & 0 \end{pmatrix}
                   \begin{pmatrix}0 & 1 \\ 1 & 0 \end{pmatrix}\otimes BA \\
                &= \begin{pmatrix}0 & 1 \\ 1 & 0 \end{pmatrix}
                   \begin{pmatrix}0 & -1 \\ 1 & 0 \end{pmatrix}\otimes (AB-BA)
                 = \begin{pmatrix}1 & 0 \\ 0 & -1 \end{pmatrix}\otimes I=A_3,\\
  A_1A_3+A_3A_1 &= \begin{pmatrix}0 & 1 \\ 1 & 0 \end{pmatrix}
                   \begin{pmatrix}1 & 0 \\ 0 & -1 \end{pmatrix}\otimes A
                 + \begin{pmatrix}1 & 0 \\ 0 & -1 \end{pmatrix} 
                   \begin{pmatrix}0 & 1 \\ 1 & 0 \end{pmatrix}\otimes A =  0, \\
  A_2A_3+A_3A_2 &= \begin{pmatrix}0 & -1 \\ 1 & 0 \end{pmatrix}
                   \begin{pmatrix}1 & 0 \\ 0 & -1 \end{pmatrix}\otimes B
                 + \begin{pmatrix}1 & 0 \\ 0 & -1 \end{pmatrix}
                   \begin{pmatrix}0 & -1 \\ 1 & 0 \end{pmatrix}\otimes B =  0. 
\end{align*}
Introducing the Pauli spin matrices
\begin{equation*}
  \sigma_1=\begin{pmatrix}0 & 1 \\ 1 & 0 \end{pmatrix}, \quad
  \sigma_2=\begin{pmatrix}0 & -i \\ i & 0 \end{pmatrix}, \quad
  \sigma_3=\begin{pmatrix}1 & 0 \\ 0 & -1 \end{pmatrix}, 
\end{equation*}
we can write $A_1=\sigma_1\otimes A$, $A_2=-i\sigma_2\otimes B$, 
and $A_3=\sigma_3\otimes I$.
Among familiar simple properties of the Pauli matrices, we have that 
 $\sigma_1\sigma_2=i\sigma_3$ and any two different Pauli matrices anticommute.
Using these relations, it follows immediately that $A_1$, $A_2$ and $A_3$ 
must satisfy (1).
\end{example}

\begin{example}
\label{exa:A_2=TA^s}
Let $s$ be a positive odd integer and take $\varphi(A)=0$ and $\psi(A)=A^s$ in
the general solution given by Theorem \ref{thm:3rel}, yielding 
\begin{equation*}
  A_1=A,\quad  A_2=T(A,B)A^s+BT(A,B),\quad A_3=T(A,B).
\end{equation*}
Considering as in Example \ref{exa:A_3=0} the case when $A=\partial_x$
and $B=M_x$ defined on the linear space ${\mathbb C}[x]$, we have here
\begin{align*}
  A_1 &= \partial_x, \\
  A_2 &= T(\partial_x,M_x)\partial_x^s+M_xT(\partial_x,M_x) = \sum_{k=0}^{\infty}
               \frac{(-2)^k}{k!}(M_x^k\partial_x^{k+s}+M_x^{k+1}\partial_x^k), \\
  A_3 &= T(\partial_x,M_x) = \sum_{k=0}^{\infty}
            \frac{(-2)^k}{k!}M_x^k\partial_x^{k}, 
\end{align*}
satisfying (\ref{rel:colLHeis}) on ${\mathbb C}[x]$ by Theorem \ref{thm:3rel}.
In a similar way as in Example \ref{exa:A_3=T}, we can define $A_1$, $A_2$ and $A_3$ 
on the space $D^{s+1}(\mathbb R)$ of all $s+1$ times differentiable functions on the 
real line, by the equations
\begin{equation*}
  (A_1f)(x) = f'(x), \quad
  (A_2f)(x) = f^{(s)}(-x)+xf(-x), \quad
  (A_3f)(x) = f(-x).
\end{equation*}
By direct computation we have
\begin{align*}
  A_1A_2f(x) &=-f^{(s+1)}(-x)+f(-x)-xf'(-x), \\
  A_2A_1f(x) &=f^{(s+1)}(-x)+xf'(-x), \\
  (A_1A_2+A_2A_1)f(x) &= f(-x)=A_3f(x), \\
  (A_1A_3+A_3A_1)f(x) &= -f'(-x)+f'(-x)=0, \\
  A_2A_3f(x) &= (-1)^sf^{(s)}(x)+xf(x), \\
  A_3A_2f(x) &= f^{(s)}(x)-xf(x), \\
  (A_2A_3+A_3A_2)f(x) &= (1+(-1)^s) f^{(s)}(x)=0. 
\end{align*}
Note that for even values of $s$, $s=2r,\ r\in \mathbb N $, we get 
operators $A_1$, $A_2$ and $A_3$ satisfying the commutation relations
\begin{equation*}
  A_1A_2 + A_2A_1 = A_3,  \quad 
  A_1A_3 + A_3A_1 = 0,    \quad
  A_2A_3 + A_3A_2 = 2A_1^{2r}.
\end{equation*}
These relations do not correspond to any color Lie algebra, though.

As in Example \ref{exa:dirsum}, let $C_0^{\infty}(\mathbb R)$ and 
$C_1^{\infty}(\mathbb R)$ be the 
subspaces of $C^{\infty}(\mathbb R)$ consisting of even and odd infinitely
differentiable functions respectively. Then we have 
 $C^{\infty}(\mathbb R)=C_0^{\infty}({\mathbb R}) \oplus C_1^{\infty}({\mathbb R})$,
and 
\begin{equation*}
  A_2:C_0^{\infty}({\mathbb R}) \rightarrow C_1^{\infty}({\mathbb R}), \qquad 
  A_2:C_1^{\infty}({\mathbb R}) \rightarrow C_0^{\infty}({\mathbb R}),
\end{equation*}
where $A_2=M_x-\partial_x^s$ on $C^{\infty}_0({\mathbb R})$ and 
$A_2=\partial_x^s-M_x$ on $C^{\infty}_1({\mathbb R})$.
Hence, we can define
\begin{equation*}
  \mathbf{A}_i:C_0^{\infty}({\mathbb R}) \oplus C_1^{\infty}({\mathbb R}) \rightarrow 
               C_0^{\infty}({\mathbb R}) \oplus C_1^{\infty}({\mathbb R}), \ i=1,2,3
\end{equation*}
where the operators $\mathbf{A}_1$, $\mathbf{A}_2$ and $\mathbf{A}_3$ are given
by the operator matrices
\begin{equation*}
  \mathbf{A}_1=\begin{pmatrix}0 & \partial_x \\ \partial_x & 0 \end{pmatrix}, \quad
  \mathbf{A}_2=\begin{pmatrix}0 &\partial_x^s -M_x \\ M_x-\partial_x^s & 0 
  \end{pmatrix}, \quad
  \mathbf{A}_3=\begin{pmatrix}I & 0 \\ 0 & -I \end{pmatrix}.
\end{equation*}
Since, for any non-negative $s$ we have 
$\partial_x(M_x-\partial_x^s)-(M_x-\partial_x^s)\partial_x=I$, 
it follows from the general result in Example \ref{exa:tensor} that the operators 
$\mathbf{A}_1$, $\mathbf{A}_2$ and $\mathbf{A}_3$ satisfy the relations 
\begin{equation*}
  \mathbf{A}_1\mathbf{A}_2+\mathbf{A}_2\mathbf{A}_1 = \mathbf{A}_3, \quad 
  \mathbf{A}_1\mathbf{A}_3+\mathbf{A}_3\mathbf{A}_1 = \mathbf{0}, \quad 
  \mathbf{A}_2\mathbf{A}_3+\mathbf{A}_3\mathbf{A}_2 = \mathbf{0},
\end{equation*} 
where $\mathbf{0}$ denotes the $2\times2$ zero matrix. 
Using the notation introduced in Example \ref{exa:dirsum} and above, 
we have that $\mathbf{A}_1$, $\mathbf{A}_2$ and $\mathbf{A}_3$ can be 
defined on $D({\mathbb R})=D_0({\mathbb R}) \oplus D_1({\mathbb R})$,
$D^s({\mathbb R})=D_0^s({\mathbb R}) \oplus D_1^s({\mathbb R})$ and
${\mathbb C}^{\mathbb R}={\mathbb C}^{\mathbb R}_0 \oplus {\mathbb C}^{\mathbb R}_1$  
respectively. For the three relations to be satisfied, it suffices to take 
$D^{s+1}({\mathbb R})=D_0^{s+1}({\mathbb R}) \oplus D_1^{s+1}({\mathbb R})$ 
as the domain of definition for all three operators.
\end{example}

\begin{example}
\label{exa:A_3=TeA-1}
Let $\varphi(A)=A$ and $\psi(A)=0$ in the general solution given by 
Theorem \ref{thm:3rel}, yielding 
\begin{equation*}
  A_1=A,\quad  A_2=BT(A,B)e^A-\tfrac{1}{2}\,T(A,B)E(A),\quad A_3=T(A,B)e^A.
\end{equation*}
Considering once again as in 
Example \ref{exa:A_3=0} the operators $A=\partial_x$
and $B=M_x$ defined on the linear space ${\mathbb C}[x]$, we have
\begin{align*}
  A_1 &= \partial_x, \\
  A_2 &= M_xT(\partial_x,M_x)e^{\partial_x}-\tfrac{1}{2}\,T(\partial_x,M_x)E(\partial_x), \\
  A_3 &= T(\partial_x,M_x)e^{\partial_x}, 
\end{align*}
or explicitly in the form of operator power series
\begin{align*}
  A_1 &= \partial_x, \\
  A_2 &=\sum_{k=0}^{\infty}\frac{(-2)^k}{k!}M_x^{k+1}\partial_x^k
        \sum_{n=0}^{\infty}\frac{\partial_x^n}{n!}-
        \tfrac{1}{2}\sum_{k=0}^{\infty}\frac{(-2)^k}{k!}M_x^k\partial_x^k
        \sum_{n=0}^{\infty}\frac{E_{2n}}{(2n)!}\partial_x^{2n}, \\
  A_3 &= \sum_{k=0}^{\infty}
         \frac{(-2)^k}{k!}M_x^k\partial_x^{k}\sum_{n=0}^{\infty}\frac{\partial_x^n}{n!}. 
\end{align*}
These operators are defined on the whole of the polynomial space
${\mathbb C}[x]$, and by Theorem \ref{thm:3rel} they satisfy 
(\ref{rel:colLHeis}) on ${\mathbb C}[x]$. 
For any differentiable function $f$, we can now define 
$A_1$, $A_2$ and $A_3$ by the equations
\begin{align*}
  (A_1f)(x) &= f'(x), \\
  (A_2f)(x) &= xf(1-x)-\tfrac{1}{2}\sum_{n=0}^{\infty}\frac{E_{2n}}{(2n)!}f^{(2n)}(-x), \\
  (A_3f)(x) &= f(1-x).
\end{align*}

By direct computation we have
\begin{multline*}
  (A_1A_2+A_2A_1)f(x) = f(1-x)-xf'(1-x)+\tfrac{1}{2}
       \sum_{n=0}^{\infty}\frac{E_{2n}}{(2n)!}f^{(2n+1)}(-x) \\
\begin{aligned}
       & +xf'(1-x)-\tfrac{1}{2}\sum_{n=0}^{\infty}
       \frac{E_{2n}}{(2n)!}f^{(2n+1)}(-x)
       =f(1-x)=A_3f(x). \\
\end{aligned}
\end{multline*}
\begin{align*}
  (A_1A_3+A_3A_1)f(x) &= -f'(1-x)+f'(1-x)=0, 
\end{align*}
\begin{multline*}
  (A_2A_3+A_3A_2)f(x) =
         xf(1-(1-x))-\tfrac{1}{2}\sum_{n=0}^{\infty}\frac{E_{2n}}{(2n)!}f^{(2n)}(1+x) \\
\begin{aligned}
           &+(1-x)f(1-(1-x))-\tfrac{1}{2}\sum_{n=0}^{\infty}\frac{E_{2n}}{(2n)!}f^{(2n)}(-(1-x)) \\
           &= f(x)-\tfrac{1}{2}\sum_{n=0}^{\infty}\frac{E_{2n}}{(2n)!}[f^{(2n)}(x-1)+f^{(2n)}(x+1)]. 
\end{aligned}
\end{multline*}
Thus, the relation $(A_2A_3+A_3A_2)f(x)=0$ is satisfied if and only if the
function $f$ satisfies
\begin{equation}
  f(x)=\tfrac{1}{2}\sum_{n=0}^{\infty}
       \frac{E_{2n}}{(2n)!}[f^{(2n)}(x-1)+f^{(2n)}(x+1)]. \label{eq:Eul-intpol}
\end{equation}
The relations (\ref{rel:colLHeis}) hold on ${\mathbb C}[x]$ and so 
(\ref{eq:Eul-intpol}) holds for $f\in{\mathbb C}[x]$.
We now give an independent combinatorial proof, showing that this formula holds 
when $f(x)$ is a polynomial in ${\mathbb C[x]}$. 
Let $f(x)=x^n$, where $n$ is a non-negative integer, and consider the sum
\begin{align*}
  s_n &= \sum_{k=0}^{\infty}\frac{E_{2k}}{(2k)!}[f^{(2k)}(x-1)+f^{(2k)}(x+1)] \\
      &=\sum_{k=0}^{\left[\frac{n}{2}\right]}E_{2k}\binom{n}{2k}[\,(x-1)^{n-2k}+(x+1)^{n-2k}\,] \\
      &=\sum_{k=0}^{\left[\frac{n}{2}\right]}E_{2k}\binom{n}{2k}
        \left[\sum_{\nu=0}^{n-2k}\binom{n-2k}{\nu}x^\nu(-1)^{n-2k-\nu}
        +\sum_{\nu=0}^{n-2k}\binom{n-2k}{\nu}x^\nu\right] \\
      &=\sum_{k=0}^{\left[\frac{n}{2}\right]}E_{2k}\binom{n}{2k}
        \sum_{\nu=0}^{n-2k}\binom{n-2k}{\nu}\left[\,1+(-1)^{n-2k-\nu}\,\right]x^\nu.
\end{align*}
Since the Euler numbers $E_{2k+1}=0$ for all non-negative integer values $k$
(see Appendix B), we can write this as
\begin{align*}
  s_n &=\sum_{m=0}^nE_m\binom{n}{m}
        \sum_{\nu=0}^{n-m}\binom{n-m}{\nu}\left[\,1+(-1)^{n-m-\nu}\,\right]x^\nu\\
      &= \sum_{\nu=0}^n\sum_{m=0}^{n-\nu}E_m\binom{n}{m}
         \binom{n-m}{\nu}\left[\,1+(-1)^{n-m-\nu}\,\right]x^\nu
\end{align*}
We have $s_n$ expressed as a polynomial of degree $n$ and with
coefficients $p_{n\nu}$ given by
\begin{align*}
  p_{n\nu} &= \sum_{m=0}^{n-\nu}E_m\binom{n}{m}
         \binom{n-m}{\nu}\left[\,1+(-1)^{n-m-\nu}\,\right]\\
            &= \frac{n!}{\nu!}\sum_{m=0}^{n-\nu}\frac{E_m}{m!(n-\nu-m)!}
               \left[\,1+(-1)^{n-\nu-m}\,\right] \\
            &=
            \frac{n!}{\nu!}\sum_{k=0}^{\left[\frac{n-\nu}{2}\right]}
             \frac{E_{2k}}{(2k)!\,(n-\nu-2k)!}\left[\,1+(-1)^{n-\nu-2k}\,\right].
\end{align*}
For odd values of $n-\nu$, we observe that $p_{n\nu}=0$ because of the
factor $1+(-1)^{n-\nu-2k}$ being zero, while for even $n-\nu$ we can
write $n-\nu=2r$ for $r\in \mathbb{N}$, and hence
\begin{align*}
  p_{n\nu} &= 2\,\frac{n!}{\nu!}\sum_{k=0}^r
              \frac{E_{2k}}{(2k)!\,(2r-2k)!}=
              \frac{2\,n!}{(2r)!\,\nu!}\sum_{k=0}^rE_{2k}\binom{2r}{2k}\\
           &= \frac{2\,n!}{(2r)!\,\nu!}\left(\frac{2^{2r+1}}{2r+1}
              (2^{2r+1}-1)B_{2r+1}+2\delta_{0,2r}\right),
\end{align*}
where $B_n$ are the Bernoulli numbers and $\delta_{k,l}$ is Kronecker's
delta. Here we have used formula (5.1.3.2) on page 385 in
\cite{PBM-3}. Some useful facts about the Bernoulli numbers can also 
be found in Appendix B.
For $\nu=n$ we have $r=0$ and
\begin{equation*}
  p_{nn}=2[2(2-1)B_1+2\delta_{00}]=4(B_1+\delta_{00})=4(-\tfrac{1}{2}+1)=2.
\end{equation*}
Since $B_{2r+1}=0$ for all $r\geq 1$, we have shown that $p_{n\nu}=0$
when $\nu<n$, and hence
\begin{equation*}
  s_n = \sum_{\nu=0}^np_{n\nu}x^{\nu}= p_{nn}x^n=2x^n. 
\end{equation*}
For every non-negative integer $n$, it holds that
\begin{equation*}
  x^n = \tfrac{1}{2}\sum_{k=0}^{\left[\frac{n}{2}\right]}
        \frac{E_{2k}}{(2k)!}\left(\frac{n!}{(n-2k)!}\,(x-1)^{n-2k}
       +\frac{n!}{(n-2k)!}\,(x+1)^{n-2k}\right).
\end{equation*}
It follows, that for every polynomial $p(x)\in \mathbb C[x]$, we obtain
\begin{equation*}
  p(x)=\tfrac{1}{2}\sum_{k=0}^{\infty}\frac{E_{2k}}{(2k)!}[p^{(2k)}(x-1)+p^{(2k)}(x+1)].
\end{equation*}

Another way to deduce this is to proceed by operator methods. 
Let us define the shift operator $S$ and the differentiation
operator $D$ as follows: For any function $f$ defined on the
real line $\mathbb R$, the action of $S$ on $f$ is given by
$(Sf)(x)=f(x+1)$. The inverse $S^{-1}$ exists and acts as
$(S^{-1}f)(x)=f(x-1)$. The operator $D$ is defined on 
$D(\mathbb R)$, the class of all differentiable functions 
on the real line, as $Df=f'$. 
If we restrict the domain of definition of both $S$ and $D$ to
the space of analytic functions on $\mathbb R$ then, as proved
in (\ref{eq:expder}), we have the equality $S=\exp(D)$. 
Moreover, it holds that
\begin{equation*}
  I+e^{2D}=I+S^2 = (S^{-1}+S)S=(S^{-1}+S)e^D. 
\end{equation*}
Using the exponential generating function of the
Euler numbers and restricting the domain of definition further
to analytic functions $f$ on $\mathbb R$, such that
\begin{equation*}
  2e^D(I+e^{2D})^{-1}f= \sum_{k=0}^{\infty}\frac{E_{2k}}{(2k)!}D^{2k}f
\end{equation*}
holds (and in particular all parts of the equality exist), then we have 
\begin{align*}
  I &=(S^{-1}+S)e^D(I+e^{2D})^{-1}
     =\tfrac{1}{2}(S^{-1}+S)\left( 2e^D(I+e^{2D})^{-1}\right) \\
    &=\tfrac{1}{2}(S^{-1}+S)\sum_{k=0}^{\infty}\frac{E_{2k}}{(2k)!}D^{2k}
     =\tfrac{1}{2}\sum_{k=0}^{\infty}\frac{E_{2k}}{(2k)!}(S^{-1}+S)D^{2k}.
\end{align*}
We have obtained
\begin{equation*}
  I =\tfrac{1}{2}\sum_{k=0}^{\infty}\frac{E_{2k}}{(2k)!}(S^{-1}D^{2k}+SD^{2k}),
\end{equation*}
or by acting on a polynomial function $p(x)$
\begin{equation*}
  p(x) =\frac{1}{2}\sum_{k=0}^{\infty}\frac{E_{2k}}{(2k)!}[p^{(2k)}(x-1)+p^{(2k)}(x+1)].
\end{equation*}
The interpolation formula (\ref{eq:Eul-intpol}) can be shown to hold
for a larger class of functions than just ${\mathbb C[x]}$, but not for all
analytic functions on $\mathbb R$. An illustrative 
example is the exponential function $e^{ax}$, where the formula can be
shown to hold for $\left|a\right|<\pi/2$ (cf. Appendix B). 
In order to extend the domain of definition for $A_1$, $A_2$ and $A_3$ from 
${\mathbb C[x]}$ to a bigger space, it would be interesting to 
characterize this class of functions and to see how the operator $A_2$ 
can be defined on a larger domain than is done here.
\end{example}

\begin{example}
\label{exa:A_3=cTeA-1}
This example contains Examples \ref{exa:A_3=0}, \ref{exa:A_3=T},
\ref{exa:A_2=TA^s} and \ref{exa:A_3=TeA-1} for specific values of the 
parameters defining $\varphi$ and $\psi$.
Let s be a positive odd integer and take $\varphi(A)=\alpha A$ 
and $\psi(A)=\beta_sA^s$ in the general solution given by 
Theorem \ref{thm:3rel}. Then we obtain 
\begin{align*}
  A_1 &= A, \quad A_3=c\,T(A,B)e^{\alpha A}, \\
  A_2 &= c\,T(A,B)E(\alpha A)[e^{\alpha A}\beta_sA^s-\tfrac{1}{2}\alpha]
       + c\,BT(A,B)e^{\alpha A}
\end{align*}
Considering once again as in 
Example \ref{exa:A_3=0} the operators $A=\partial_x$
and $B=M_x$ defined on the linear space ${\mathbb C}[x]$, we have
\begin{align*}
  A_1 &= \partial_x, \\
  A_2 &= c\,T(\partial_x,M_x)E(\alpha \partial_x)[e^{\alpha \partial_x}\beta_s\partial_x^s-\tfrac{1}{2}\alpha] + c\,M_xT(\partial_x,M_x)e^{\alpha \partial_x}, \\
  A_3 &= c\,T(\partial_x,M_x)e^{\alpha \partial_x}. 
\end{align*}

We can now define $A_1$, $A_2$ and $A_3$ on the polynomial space 
${\mathbb C}[x]$ by the equations
\begin{align*}
  (A_1f)(x) &= f'(x), \\
  (A_2f)(x) &= c\,T(\partial_x,M_x)E(\alpha \partial_x)[\beta_sf^{(s)}(x+\alpha)
             -\tfrac{1}{2}\alpha f(x)] + cxf(-x+\alpha) \\
            &= c\,T(\partial_x,M_x)\sum_{n=0}^{\infty}\frac{E_{2n}}{(2n)!}\alpha^{2n}
               [\beta_sf^{(s+2n)}(x+\alpha)-\tfrac{1}{2}\alpha f^{(2n)}(x)] \\
            & \qquad + cxf(-x+\alpha) \\
            &= c\sum_{n=0}^{\infty}\frac{E_{2n}}{(2n)!}\alpha^{2n}
               [\beta_sf^{(s+2n)}(-x+\alpha)-\tfrac{1}{2}\alpha f^{(2n)}(-x)] \\
            & \qquad + cxf(-x+\alpha) \\
  (A_3f)(x) &= c\,T(\partial_x,M_x)e^{\alpha \partial_x}f(x) = cf(-x+\alpha).
\end{align*}

By direct computation we have
\begin{align*}
  (A_1A_2f)(x) &= c\sum_{n=0}^{\infty}\frac{E_{2n}}{(2n)!}\alpha^{2n}
               [-\beta_sf^{(s+2n+1)}(-x+\alpha)+\tfrac{1}{2}\alpha f^{(2n+1)}(-x)] \\
               & \qquad -  cxf'(-x+\alpha)+ cf(-x+\alpha).
\end{align*}

\begin{align*}
  (A_2A_1f)(x) &= c\sum_{n=0}^{\infty}\frac{E_{2n}}{(2n)!}\alpha^{2n}
               [\beta_sf^{(s+2n+1)}(-x+\alpha)-\tfrac{1}{2}\alpha f^{(2n+1)}(-x)] \\
               & \qquad +  cxf'(-x+\alpha).
\end{align*}

\begin{equation*}
  (A_1A_2+A_2A_1)f(x) = cf(-x+\alpha)=(A_3f)(x), 
\end{equation*}

\begin{equation*}
  (A_1A_3+A_3A_1)f(x) = -cf'(-x+\alpha)+cf'(-x+\alpha)=0, 
\end{equation*}

\begin{align*}
  (A_2A_3f)(x) &= c\sum_{n=0}^{\infty}\frac{E_{2n}}{(2n)!}\alpha^{2n}
               [c\beta_s(-1)^sf^{(s+2n)}(x)-\tfrac{1}{2}c\alpha f^{(2n)}(x-\alpha)] \\
               & \qquad +  c^2xf(x).
\end{align*}

\begin{align*}
  (A_3A_2f)(x) &= c^2\sum_{n=0}^{\infty}\frac{E_{2n}}{(2n)!}\alpha^{2n}
               [\beta_sf^{(s+2n)}(x)-\tfrac{1}{2}\alpha f^{(2n)}(x+\alpha)] \\
               & \qquad +  c^2(-x+\alpha)f(x).
\end{align*}

\begin{multline*}
  (A_2A_3+A_3A_2)f(x) =
                       c^2xf(x)-\tfrac{1}{2}c^2\sum_{n=0}^{\infty}\frac{E_{2n}}{(2n)!}\alpha^{2n+1}f^{(2n)}(x-\alpha) \\
\begin{aligned}
           &+c^2(-x+\alpha)f(x)-\tfrac{1}{2}c^2\sum_{n=0}^{\infty}\frac{E_{2n}}{(2n)!}\alpha^{2n+1}f^{(2n)}(x+\alpha) \\
           &= \alpha c^2 f(x)-\tfrac{1}{2}c^2\sum_{n=0}^{\infty}\frac{E_{2n}}{(2n)!}\alpha^{2n+1}[f^{(2n)}(x-\alpha)+f^{(2n)}(x+\alpha)]. 
\end{aligned}
\end{multline*}
Thus, the relation $(A_2A_3+A_3A_2)f(x)=0$ is satisfied if and only if the
function $f$ satisfies
\begin{equation}
  f(x)=\tfrac{1}{2}\sum_{n=0}^{\infty}
       \frac{E_{2n}}{(2n)!}\alpha^{2n}[f^{(2n)}(x-\alpha)+f^{(2n)}(x+\alpha)]. \label{eq:Eul-intpolalf}
\end{equation}
\end{example}

\begin{example}
\label{exa:A_3=TeA-2}
Let us have a closer look at the operators studied in 
Example \ref{exa:A_3=TeA-1}, namely
\begin{align*}
  A_1 &= \partial_x, \\
  A_2 &= M_xT(\partial_x,M_x)e^{\partial_x}-\tfrac{1}{2}\,T(\partial_x,M_x)E(\partial_x), \\
  A_3 &= T(\partial_x,M_x)e^{\partial_x}, 
\end{align*}
We can express $E(\partial_x)$ in terms of the generating function for
the Euler numbers, obtaining
\begin{align*}
  E(\partial_x) &= 2e^{\partial_x}(e^{2\partial_x}+1)^{-1} 
                 = e^{\partial_x}(\tfrac{1}{2}+\tfrac{1}{2}e^{2\partial_x})^{-1}                    
                 = e^{\partial_x}(1+\tfrac{1}{2}(e^{2\partial_x}-1))^{-1} \\
                &= e^{\partial_x}
                   \sum_{k=0}^{\infty}\frac{(-1)^k}{2^k}(e^{2\partial_x}-1)^k 
                 = e^{\partial_x}\sum_{k=0}^{\infty}\frac{(-1)^k}{2^k}
                   \sum_{l=0}^k\binom{k}{l}e^{2(k-l)\partial_x}(-1)^l \\
                &= \sum_{k=0}^{\infty}\frac{(-1)^k}{2^k}
                   \sum_{l=0}^k(-1)^l\binom{k}{l}e^{(2(k-l)+1)\partial_x} 
\end{align*}
By virtue of equations (\ref{eq:expder}) and (\ref{eq:parity}) in 
Example \ref{exa:A_3=TeA-1}, it is now reasonable to
define 
\begin{align*}
  (A_1f)(x) &= f'(x), \\
  (A_2f)(x) &= xf(1-x)-\sum_{k=0}^{\infty}\frac{(-1)^k}{2^{k+1}}
              \sum_{l=0}^k(-1)^l\binom{k}{l}f(2(k-l)+1-x), \\
  (A_3f)(x) &= f(1-x),
\end{align*}
where $f$ is a polynomial ${\mathbb C} [x]$.
In order to verify the commutation relations, we calculate
\begin{multline*}
  A_1A_2f(x) = \partial_x \left(xf(1-x)-\sum_{k=0}^{\infty}
               \sum_{l=0}^k\frac{(-1)^{k+l}}{2^{k+1}}\binom{k}{l}f(2(k-l)+1-x)\right) \\
\begin{aligned}   
             &= f(1-x)-xf'(1-x)+\sum_{k=0}^{\infty}\sum_{l=0}^k
             \frac{(-1)^{k+l}}{2^{k+1}}\binom{k}{l}f'(2(k-l)+1-x)
\end{aligned}
\end{multline*}
\begin{equation*}
  A_2A_1f(x) =xf'(1-x)+\sum_{k=0}^{\infty}\sum_{l=0}^k
               \frac{(-1)^{k+l}}{2^{k+1}}\binom{k}{l}f'(2(k-l)+1-x)
\end{equation*}
We have then
\begin{align*}
  (A_1A_2+A_2A_1)f(x) &= f(1-x)=A_3f(x), \\
  (A_1A_3+A_3A_1)f(x) &= -f'(1-x)+f'(1-x)=0.  
\end{align*}
Moreover
\begin{align*}
  A_2A_3f(x) &= A_2f(1-x)=         
             xf(x)-\sum_{k=0}^{\infty}\sum_{l=0}^k
                \frac{(-1)^{k+l}}{2^{k+1}}\binom{k}{l}f(x-2(k-l)) \\
  A_3A_2f(x) &=A_3\left(xf(1-x)-\sum_{k=0}^{\infty}\sum_{l=0}^k
               \frac{(-1)^{k+l}}{2^{k+1}}\binom{k}{l}f(2(k-l)+1-x)\right) \\
             &= (1-x)f(x)-\sum_{k=0}^{\infty}\sum_{l=0}^k
                \frac{(-1)^{k+l}}{2^{k+1}}\binom{k}{l}f(x+2(k-l))
\end{align*}
and hence
\begin{multline*}
  (A_2A_3+A_3A_2)f(x) =xf(x)+(1-x)f(x) \\
\begin{aligned}
              &-\sum_{k=0}^{\infty}\sum_{l=0}^k
                \frac{(-1)^{k+l}}{2^{k+1}}\binom{k}{l}
                [f(x-2(k-l))+f(x+2(k-l))] \\
              &= f(x)-\sum_{k=0}^{\infty}\sum_{l=0}^k
                 \frac{(-1)^{k+l}}{2^{k+1}}\binom{k}{l}
                 [f(x-2(k-l))+f(x+2(k-l))]
\end{aligned}
\end{multline*}
In order to satisfy the relation $(A_2A_3+A_3A_2)f(x)=0$, we must have
\begin{equation}
  f(x)=\sum_{k=0}^{\infty}\frac{(-1)^k}{2^{k+1}} 
       \sum_{l=0}^k(-1)^l\binom{k}{l}[f(x-2(k-l))+f(x+2(k-l))]. \label{eq:diffintpol}
\end{equation}
We shall now demonstrate that this formula holds if $f(x)$ is an
arbitrary polynomial in $x$. For this purpose we need to use some
properties of the Stirling numbers of the second kind, see Chapter 5
in \cite{Cam}. We take 
$f(x)=x^n$ with $n$ a non-negative integer and consider 
the sum
\begin{align*}
  \sigma_n &= \sum_{k=0}^{\infty}\frac{(-1)^k}{2^k}
              \sum_{l=0}^k(-1)^l\binom{k}{l}[f(x-2(k-l))+f(x+2(k-l))] \\
           &= \sum_{k=0}^{\infty}\frac{(-1)^k}{2^k}
              \sum_{l=0}^k(-1)^l\binom{k}{l}[(x-2(k-l))^n+(x+2(k-l))^n] \\
           &= \sum_{k=0}^{\infty}\frac{(-1)^k}{2^k}
              \sum_{l=0}^k(-1)^l\binom{k}{l}
              \sum_{\nu=0}^n\binom{n}{\nu}x^{n-\nu}2^\nu(k-l)^\nu(1+(-1)^\nu) \\
           &= \sum_{\nu=0}^n\binom{n}{\nu}2^\nu(1+(-1)^\nu)
              \sum_{k=0}^{\infty}\frac{(-1)^k}{2^k}\sum_{l=0}^k
              (-1)^l\binom{k}{l}(k-l)^\nu x^{n-\nu} \\
           &= \sum_{\nu=0}^n\binom{n}{\nu}2^\nu(1+(-1)^\nu)
              \sum_{k=0}^{\infty}\frac{(-1)^k}{2^k}\sum_{l=0}^k
              (-1)^{k-l}\binom{k}{l}l^\nu x^{n-\nu} \\
           &= \sum_{\nu=0}^n\binom{n}{\nu}2^\nu(1+(-1)^\nu)
              \sum_{k=0}^{\infty}\frac{(-1)^kk!}{2^k}
              S(\nu,k)x^{n-\nu},
\end{align*}
where we have introduced the Stirling numbers of the second kind,
given by
\begin{align*}
  S(m,k) &= \frac{1}{k!}\sum_{j=0}^k(-1)^{k-j}\binom{k}{j}j^m, 
           \quad m\geq 1,\  k\geq 0, \\
  S(0,0) &= 1, \quad S(0,k)=0, \quad k\geq 1.     
\end{align*}
The sum $\sigma_n$ is now expressed as a polynomial of degree
$n$ with coefficients $q_{n\nu}$ given as
\begin{equation*}
  q_{n\nu} = \sum_{\nu=0}^n\binom{n}{\nu}2^\nu(1+(-1)^\nu)
              \sum_{k=0}^\nu\frac{(-1)^kk!}{2^k}S(\nu,k)     
\end{equation*} 
since for $k>m$, we have by equation (4.2.2.3) on page 608 in
\cite{PBM-1}
\begin{equation*}
  \sum_{j=0}^k(-1)^{k-j}\binom{k}{j}j^m=0.     
\end{equation*} 
For odd values of $\nu$, we note that $q_{n\nu}=0$ due to the 
factor $(1+(-1)^\nu)$. Now define $f_0=g_0=0$ and let for all
positive integers $k,m$
\begin{equation*}
  f_k = \frac{(-1)^kk!}{2^k}, \qquad 
  g_m = \sum_{k=0}^mS(m,k)f_k.    
\end{equation*}
The exponential generating functions $F(t)$ and $G(t)$ corresponding
to the sequences $(f_k)$ and $(g_m)$ respectively, are defined as
\begin{equation*}
  F(t) = \sum_{k=0}^{\infty}f_k\frac{t^k}{k!}, \qquad
  G(t) = \sum_{m=0}^{\infty}g_m\frac{t^k}{m!}.
\end{equation*}
We find
\begin{equation*}
  F(t) = \sum_{k=0}^{\infty}f_k\frac{t^k}{k!}=
         \sum_{k=1}^{\infty}\frac{(-1)^k}{2^k}t^k=
         \sum_{k=1}^{\infty}\left(\frac{-t}{2}\right)^k
         =\frac{-t}{2+t}. 
\end{equation*}
By Theorem 5.4.2 in \cite{Cam} it follows that
\begin{equation*}
  G(t) = F(e^t-1)=\frac{1-e^t}{1+e^t}, 
\end{equation*}
and since
\begin{equation*}
  G(-t) =\frac{1-e^{-t}}{1+e^{-t}}=\frac{e^t-1}{e^t+1}=-G(t), 
\end{equation*}
$G(t)$ is an odd generating function, and hence for positive even
values of $\nu$ 
\begin{equation*}
  g_{\nu} =  \sum_{k=0}^{\nu}S(\nu,k)f_k=
          \sum_{k=0}^\nu\frac{(-1)^kk!}{2^k}S(\nu,k)=0, 
\end{equation*}
yielding the result
\begin{equation*}
  \sigma_n = \sum_{\nu=0}^nq_{n\nu}x^{n-\nu}= q_{n0}x^n=2x^n. 
\end{equation*}
This proves that, for any non-negative integer $n$, we have
\begin{equation*}
  x^n = \sum_{k=0}^{\infty}\frac{(-1)^k}{2^{k+1}}
              \sum_{l=0}^k(-1)^l\binom{k}{l}[(x-2(k-l))^n+(x+2(k-l))^n].
\end{equation*}
It follows, that for every polynomial $q(x)\in \mathbb C[x]$, we obtain
\begin{equation*}
  q(x)=\sum_{k=0}^{\infty}\frac{(-1)^k}{2^{k+1}} 
       \sum_{l=0}^k(-1)^l\binom{k}{l}[q(x-2(k-l))+q(x+2(k-l))].
\end{equation*}
The formula (\ref{eq:diffintpol}) holds for other functions than 
polynomials. An interesting task would be to find a characterization 
of this class of functions.
\end{example}

\begin{acknowledgements}
\markboth{Acknowledgements}{}
The authors would like to thank Professor Bernt Lindstr{\"o}m for 
suggesting the main ideas behind the combinatorial
argument in Example \ref{exa:A_3=TeA-2}.
Thanks are also due to Dr.~Lars Hellstr{\"o}m and 
Dr.~Edwin Langmann for useful comments.
The first author is grateful to the Department of Mathematics
at the Lund Institute of Technology for hospitality during his
visits in Lund. 
\\
\\
\end{acknowledgements}
\newpage

\section*{Appendix A}
\label{sec:appA}
\markboth{Appendix B}{Appendix A}
Recall that a \mbox{$\mathbb{Z}_{2}^{n}$-graded} (color)
generalized Lie algebra is a \mbox{$\mathbb{Z}_{2}^{n}$-graded} linear
space $$X=\bigoplus_{\gamma \in \mathbb{Z}_{2}^{n}}
{X_{\gamma}}$$ with a bilinear multiplication 
(bracket)
\mbox{$\langle\,\cdot \, , \cdot\,\rangle: X \times X \to X$} 
obeying:

\vbox{\begin{description}
 \item[Grading axiom: ] 
$\langle\, X_{\alpha}\, , X_{\beta}\,\rangle \subseteq  X_{\alpha +
\beta}$. 
 \item[Graded skew-symmetry:] 
$\langle\, a\, , b\,\rangle=- (-1)^{\alpha \cdot \beta}\langle \,b\, , a\,\rangle$. 
 \item[Generalized Jacobi identity:]
\begin{equation*}
\hspace{-1.3cm} (-1)^{\alpha \cdot \gamma}\langle\, a\, , \langle\, b\, , c\,\rangle\,\rangle \mbox{} +
(-1)^{\gamma \cdot \beta}\langle\, c\, , \langle\, a\, , b\,\rangle\,\rangle 
\mbox{} + (-1)^{
\beta \cdot \alpha}\langle\, b\, , \langle\, c\, , a\,\rangle\,\rangle=0 
\end{equation*}
\end{description}}
\noindent for all
$\alpha=(\alpha_{1},\dots,\alpha_{n})$ , 
$\beta=(\beta_{1},\dots,\beta_{n})$ ,  
$\gamma=(\gamma_{1},\dots,\gamma_{n})$ 
in $\mathbb{Z}_{2}^{n},$  and 
\mbox{$a \in X_{\alpha}$ , $b \in
X_{\beta}$ , $c \in X_{\gamma}$}, 
where $\alpha \cdot \beta =\sum_{i=1}^n \alpha_i\beta_i$ etc., with
$\sum$ meaning
 addition in 
$\mathbb{Z}_{2}$. 
The elements of  
$\bigcup_{\gamma \in \mathbb{Z}_{2}^{n}} {X_{\gamma}}$ are called
homogeneous.

Any \mbox{$\mathbb{Z}_{2}^{n}$-graded} generalized 
Lie algebra $X$ can be embedded in its universal
enveloping algebra $U(X)$ in such a way that, 
for homogeneous \mbox{$a \in X_{\alpha}$} and  
$b \in X_{\beta}$, the bracket $\langle\,\cdot \, , \cdot\,\rangle$
becomes  a commutator $[a\, , b]=ab-ba$ when 
$\alpha \cdot \beta$ 
is even, or an
anticommutator $\{a\, ,b\}=ab+ba$ when 
$\alpha \cdot \beta$
is odd \cite{Sh-GLA}.  

Now take $X$ to be a \mbox{$\mathbb{Z}_{2}^{3}$-graded} linear space
$$X = X_{(1,1,0)} \oplus X_{(1,0,1)} 
\oplus X_{(0,1,1)}$$ with the homogeneous basis 
 $A_1 \in X_{(1,1,0)}$, $A_2 \in X_{(1,0,1)}$, 
$A_3 \in X_{(0,1,1)}$. The homogeneous
components graded by the elements of
$\mathbb{Z}_{2}^{3}$ different from  $(1,1,0)$, $(1,0,1)$ and 
$(0,1,1)$  are zero and so are omitted. 
If the \mbox{$\mathbb{Z}_{2}^{3}$-graded} bilinear
multiplication $\langle\,\cdot \, , \cdot\,\rangle$ turns $X$ into
a \mbox{$\mathbb{Z}_{2}^{3}$-graded} generalized Lie algebra,
then  $ \langle\, A_{i}\, , A_{i}\,\rangle=0 , \ i = 1,2,3$  and
\begin{equation*}
\langle\, A_{1}\, , A_{2}\,\rangle=c_{12} A_3, \quad
\langle\, A_{2}\, , A_{3}\,\rangle=c_{23} A_1, \quad
\langle\, A_{3}\, , A_{1}\,\rangle=c_{31} A_2 \ .
\end{equation*} 
When $a$ and $b$ are in different homogeneous subspaces,
it follows that $\langle\, a\, , b\,\rangle=\langle\, b\, , a\,\rangle$,  
whereas $\langle\, a\, , b\,\rangle=-\langle\, b\, , a\,\rangle$
if $a$ and $b$ belong to the same one. Moreover, the
generalized Jacobi identity is valid.
Now put $c_{12}=1$, $c_{23}=0$ and $c_{31}=0$. 
The algebra $X$ so defined has as its universal
enveloping algebra the color analogue of the  Heisenberg
Lie algebra.
\newpage
\section*{Appendix B}
\label{sec:appB}

The Bernoulli polynomials $B_k(x)$ can be defined in terms of
their exponential generating function
\begin{equation*}
  \frac{te^{xt}}{e^t-1}=\sum_{k=0}^{\infty}B_k(x)\frac{t^k}{k!}\,, 
   \quad \left|t\right|<2\pi.  
\end{equation*}
The four polynomials of lowest degree are
\begin{equation*}
  B_0(x)=1, \quad B_1(x) =x-\tfrac{1}{2}, \quad B_2(x)
  =x^2-x+\tfrac{1}{6}, \quad B_3(x)=x^3-\tfrac{3}{2}x^2+\tfrac{1}{2}x.  
\end{equation*}
The Bernoulli numbers $B_k$ are then defined as the values of 
$B_k(x)$ at the origin,
$B_k=B_k(0)$, from which it follows that $B_0=1$,
$B_1=-\tfrac{1}{2}$, and for $k\geq 1$
\begin{equation*}
  B_{2k+1}=0, \quad B_{2k}=(-1)^{k+1}
           \frac{2\,(2k)!}{\pi^{2k}(2^{2k}-1)}
           \sum_{\nu=0}^{\infty}(2\nu+1)^{-2k}.  
\end{equation*}

In a similar way, one can define a sequence of polynomials
$E_k(x)$, called the Euler polynomials, by specifying their
exponential generating function as
\begin{equation*}
  \frac{2e^{xt}}{e^t+1}=\sum_{k=0}^{\infty}E_k(x)\frac{t^k}{k!}\,, 
   \quad \left|t\right|<\pi.  
\end{equation*}
The four polynomials of lowest degree are
\begin{equation*}
  E_0(x)=1, \quad E_1(x) =x-\tfrac{1}{2}, \quad E_2(x)
  =x^2-x, \quad E_3(x)=x^3-\tfrac{3}{2}x^2+\tfrac{1}{4}.  
\end{equation*}
The Euler numbers $E_k$ are then defined as the integers 
$E_k=2^kE_k(\frac{1}{2})$. It follows that $E_0=1$, $E_1=0$,
 $E_2=-1$, $E_3=0$, and generally for $k\geq 0$
 \begin{equation*}
  E_{2k+1}=0, \quad E_{2k}=(-1)^k
           \frac{(2k)!\,2^{2k+2}}{\pi^{2k+1}}
           \sum_{\nu=0}^{\infty}(-1)^\nu(2\nu+1)^{-2k-1}.  
\end{equation*}
The Euler numbers have an exponential generating function
obtained by setting $x=1/2$ and replacing $t$ by $2t$ in the 
exponential generating function of
the Euler polynomials
\begin{equation*}
 E(t)= \frac{2e^t}{e^{2t}+1}=\sum_{k=0}^{\infty}E_k\frac{t^k}{k!}.  
\end{equation*}


\newpage

\noindent
\begin{thebibliography}{99} 
\markboth{References}{References}
\bibitem{Agrav} {\rm V.K. Agrawala,}  
{\em Invariants of generalized Lie algebras,}   
{\rm Hadronic J. } {\bf 4} (1981), 
444--496.

\bibitem{BahtMPZ} {\rm Yu.A. Bahturin, A.A. Mikhalev, 
V.M. Petrogradsky, M.V. Zaicev,} 
{\em Infinite Dimensional Lie Superalgebras,}  
{\rm Walter de Gruyter, Berlin, 1992.}

\bibitem{Blo} {\rm R.E. Block,} 
{\em The Irreducible Representations of the Lie Algebra
$\mathfrak{sl}$(2) and of the Weyl Algebra,}  
{\rm Adv. Math.} {\bf 39} (1981), 69--110.

\bibitem{Cam} {\rm P.J. Cameron,} 
{\em Combinatorics: topics, techniques, algorithms},  
{\rm Cambridge University Press}, 1994.

\bibitem{CNS} {\rm L. Corwin, Y. Ne'eman, S. Sternberg,} 
{\em  Graded Lie algebras in mathematics and physics (Bose-Fermi
symmetry),}
{\rm  Rev. Modern Phys.} 
{\bf47} (1975), 573--603. 

\bibitem{Gr-Jar} {\rm H.S. Green, P.D. Jarvis,}  
{\em Casimir invariants, characteristic 
identities and Young diagrams for  
Colour algebras and superalgebras,}   
{\rm  J. Math. Phys.}  {\bf 24} (1983), 1681--1687. 

\bibitem{HelSil-book} 
{\rm L. Hellstr{\"o}m, S.D. Silvestrov,}
\emph{Commuting Elements in q-Deformed Heisenberg Algebras},
World Scientific, 2000, 256 pp.\\
 (ISBN: 981-02-4403-7).

\bibitem{Kac-art-adv} {\rm V.G. Kac,}  
{\em Lie Superalgebras,}   
{\rm  Adv. Math.} {\bf 26} (1977), 
8--96.

\bibitem{K-M-NPR} {\rm  M.V. Karasev, V.P. Maslov,} 
{ \em Non-Lie permutation relations,}  
{ \rm Russian Math. Surveys } 
{\bf45} (1990), 51--98.

\bibitem{Kle-cfac} {\rm R. Kleeman,}  
{\em Commutation factors on generalized Lie algebras,}   
{\rm  J. Math. Phys.}   {\bf 26} (1985), 
2405--2412.

\bibitem{Kwasn-2} {\rm A.K. Kwasniewski,} 
{\em  Clifford- and Grassmann-like algebras -- Old and new,}
{\rm  J. Math. Phys.} 
{\bf26} (1985), 2234--2238. 

\bibitem{Kwasn-3} {\rm A.K. Kwasniewski,}  
{\em On Graded Lie-like Algebras,}   
{\rm Bulletin de la {Soci}\'e{t}\'e 
des sciences et des lettres 
de \L{\'o}{d}\'z }  {\bf 39}, 6 (1989).

\bibitem{LS} {\rm H. Ljungqvist, S.D. Silvestrov,
{\em Involutions in three-dimensional coloured Lie algebras,}
Research Reports {\bf 6}, Dep. Math., Ume{\aa} University, (1996), 36 pp.}

\bibitem{LucRit} {\rm J. Lukierski, V. Rittenberg,}  
{\em  Color-de Sitter and color-conformal superalgebras,}   
{\rm  Phys. Rev. D}  {\bf18} (1978), 385--389.

\bibitem{Marc-colext} {\rm  W. Marcinek,}  
{\em Colour extensions of Lie algebras and superalgebras,}   
{\rm Preprint {\bf 746}, University of Wroc\l{aw}}  (1990). 

\bibitem{Marc-GLA} {\rm W. Marcinek,} 
{\em  Generalized Lie algebras and Related Topics,1,2,}
{\rm   Acta Univ. Wratislaviensis 
( Matematyka, Fizyka, Astronomia )} {\bf LV}, 1170 (1991), 3--52.

\bibitem{Mosol} {\rm  M.V. Mosolova,} 
{\em Functions of non-commuting operators that generate a
graded Lie algebra,} {\rm Mat. Zametki} 
{\bf 29} (1981), 34--45.

\bibitem{OstSam-book}
{\rm V.L. Ostrovsky\u{\i}, Yu.S. Samo\u{\i}lenko,}
{\em Introduction to the Theory of Representations of Finitely
Presented $*$-Algebras. I. Representations by bounded operators},
{\rm Rev. Math. Math. Phys. {\bf 11},
The Gordon and Breach Publ. Group, 1999.}

\bibitem{O-Sil} {\rm  V.L. Ostrovskii, S.D. Silvestrov,} 
{\em Representations of the real forms of the graded
analogue of the Lie algebra $sl(2,\mathbb{C})$,}
{\rm Ukrain.\ Mat.\ Zh.\ 
{\bf 44} (1992), 1518--1524; 
(English translation: Ukrainian Math. J. {\bf 44} (1993), 1395--1401).} 

\bibitem{PSS} {\rm L. Persson, S.D. Silvestrov, P. Strunk}, 
{\em Central elements of the second order in three-dimensional
generalised Lie algebras}, {\em Czechoslovak Journal of Physics} 
{\bf 47} (1997), 99--106. \\
{\rm (Also published as Research Report {\bf 5}, 
Dep. Math., Ume{\aa} University, (1996), \mbox{24 pp.)}} 

\bibitem{P} {\rm K.L. Price}, 
{\em 
Primeness Criteria for Universal Enveloping Algebras of Lie Color
Algebras}, {\em Journal of Algebra} 
{\bf 235} (2001), 589--607.
 
\bibitem{PBM-1} {\rm  A.P. Prudnikov, Yu.A. Brychkov, O.I. Marichev,}  
{ \em Integrals and Series, Vol. 1: Elementary Functions, }   
{ \rm Gordon and Breach Science Publishers, 1990.}  
{\rm (Translated from the Russian original edition: Nauka, Moscow, 1986.)} 

\bibitem{PBM-3} {\rm  A.P. Prudnikov, Yu.A. Brychkov, O.I. Marichev,}  
{ \em Integrals and Series, Vol. 3: More Special Functions, }   
{ \rm Gordon and Breach Science Publishers, 1990.}  
{\rm (Translated from the Russian original edition: Nauka, Moscow, 1986.)} 

\bibitem{Put} {\rm C.R. Putnam,} 
{\em  Commutation Properties of Hilbert Space Operators and Related Topics,}  
{\rm Springer-Verlag, Berlin Heidelberg}, 1967.

\bibitem{R-W(GL)} {\rm  V. Rittenberg, D. Wyler,}  
{ \em Generalized Superalgebras,}   
{ \rm Nucl. Phys. }  
{\bf B139} (1978), 189--202. 

\bibitem{Beta} {\rm L. R{\aa}de, B. Westergren,} 
{\em  Mathematics Handbook for Science and Engineering,}  
{\rm Studentlitteratur, Lund}, 1995.

\bibitem{S} {\rm Yu.S. Samoilenko,} 
{\em  Spectral Theory of Families of Self-Adjoint Operators,}  
{\rm  Kluwer Acad. Publ., Dordrecht, 1991.} 

\bibitem{Sh-GLA} {\rm M. Scheunert,}
{\em Generalized Lie algebras,}
{\rm  J. Math. Phys.} 
{\bf 20} (1979), 712--720. 

\bibitem{Sh-GTC} {\rm M. Scheunert,}  
{\em  Graded tensor calculus,}   
{\rm J. Math. Phys. }  {\bf 24} (1983), 
2658--2670.

\bibitem{Sh-CsmrGLA} {\rm M. Scheunert,}  
{\em Casimir elements of $\varepsilon$-Lie algebras,}   
{\rm J. Math. Phys. } {\bf 24} (1983), 
 2671--2680.

\bibitem{SilSig} {\rm G.~Sigurdsson, S.D.~Silvestrov,}
{\em Canonical involutions in three-dimensional generalised Lie algebras,
Czechoslovak Journal of Physics,} {\bf 50} (2000), 181--186.

\bibitem{S-PhD} {\rm S.D. Silvestrov,}
{\em Representations of Commutation Relations.
A Dynamical Systems Approach},
{\rm Hadronic Journal Supplement}, {\bf 11} (1996), 1--116.

\bibitem{S-stm} {\rm S.D. Silvestrov,}
{\em  Hilbert space representations of 
the graded analogue of the Lie algebra of the group of plane
motions,} {\rm Studia Mathematica} {\bf 117} (1996), 195--203. 

\bibitem{S-clas} {\rm S.D. Silvestrov,} 
{\em On the classification of 
\mbox{3-dimensional} coloured Lie algebras,} 
in {\em Quantum Groups and Quantum Spaces}, 
{\rm Banach Center Publications {\bf 40} (1997), 159--170.} 

\bibitem{NW-MaQuant} {\rm N. Weaver,}  
{\em Mathematical Quantization,}   
{\rm (Studies in advanced mathematics), Chapman \& Hall/CRC, 2001.}

\bibitem{Wie-Qua} {\rm H. Wielandt,}  
{\em  {\"U}ber die Unbeschr{\"a}nkheit der Schr{\"o}dingerschen
Operatoren der Quanten-mechanik,}   
{\rm Math. Ann.}  {\bf 121} (1949), 21.

\bibitem{Win-qm} {\rm A. Wintner,}  
{\em  The unboundedness of quantum-mechanical matrices,}   
{\rm Phys. Rev.}  {\bf 71} (1947), 
738--739.

\end{thebibliography}
\end{document}